\documentclass[11pt]{amsart}

\usepackage{enumerate,url,amssymb,  mathrsfs, nicefrac, pifont, stmaryrd, graphicx, pdfsync, upref, color,soul}
\newtheorem{theorem}{Theorem}[section]
\newtheorem{lemma}[theorem]{Lemma}
\newtheorem*{lemma*}{Lemma}

\newtheorem{proposition}[theorem]{Proposition}
\newtheorem{corollary}[theorem]{Corollary}

\theoremstyle{definition}
\newtheorem{definition}[theorem]{Definition}
\newtheorem{example}[theorem]{Example}

\newtheorem{question}[theorem]{Question}

\theoremstyle{remark}
\newtheorem{remark}[theorem]{Remark}

\numberwithin{equation}{section}

\newcommand{\abs}[1]{\lvert#1\rvert}
\newcommand{\Abs}[1]{\Big\lvert#1\Big\rvert}
\newcommand{\norm}[1]{\lVert#1\rVert}

\newcommand{\C}{\mathbb{C}}

\newcommand{\W}{\mathscr{W}}

\newcommand{\E}{\mathcal{E}}
\newcommand{\EE}{\mathsf{E}}

\newcommand{\R}{\mathbb{R}}
\newcommand{\X}{\mathbb{X}}
\newcommand{\U}{\mathbb{U}}

\newcommand{\s}{\mathbb{S}}

\newcommand{\Y}{\mathbb{Y}}

\newcommand{\Z}{\mathbb{Z}}
\newcommand{\T}{\mathbb{T}}

\newcommand{\MPS}{\mathcal{MPS}}
\newcommand{\Ho}{\mathscr{H}}

\newcommand{\dtext}{\textnormal d}
\newcommand{\cto}{\xrightarrow[]{{}_{cd}}}
\newcommand{\onto}{\xrightarrow[]{{}_{\!\!\textnormal{onto}\!\!}}}
\newcommand{\into}{\xrightarrow[]{{}_{\!\!\textnormal{into}\!\!}}}
\DeclareMathOperator{\diam}{diam}
\DeclareMathOperator{\esssup}{ess\, sup}
\DeclareMathOperator{\dist}{dist}

\DeclareMathOperator{\re}{Re}
\DeclareMathOperator{\im}{Im}

\DeclareMathOperator{\loc}{loc}
\DeclareMathOperator{\cd}{cd}
\DeclareMathOperator{\Mod}{Mod}
\def\Xint#1{\mathchoice
{\XXint\displaystyle\textstyle{#1}}%
{\XXint\textstyle\scriptstyle{#1}}%
{\XXint\scriptstyle\scriptscriptstyle{#1}}%
{\XXint\scriptscriptstyle\scriptscriptstyle{#1}}%
\!\int}
\def\XXint#1#2#3{{\setbox0=\hbox{$#1{#2#3}{\int}$}\vcenter{\hbox{$#2#3$}}\kern-.5\wd0}}
\def\dashint{\Xint-}
\def\Xiint#1{\mathchoice
{\XXiint\displaystyle\textstyle{#1}}%
{\XXiint\textstyle\scriptstyle{#1}}%
{\XXiint\scriptstyle\scriptscriptstyle{#1}}%
{\XXiint\scriptscriptstyle\scriptscriptstyle{#1}}%
\!\iint}
\def\XXiint#1#2#3{{\setbox0=\hbox{$#1{#2#3}{\iint}$}\vcenter{\hbox{$#2#3$}}\kern-.5\wd0}}
\def\dashiint{\Xiint{-\!-}}
\def\le{\leqslant}
\def\ge{\geqslant}
\begin{document}

\title[The Hopf-Laplace equation]{The~Hopf-Laplace~equation: harmonicity~and~regularity}

\author[Cristina]{Jan Cristina}
\address{University of Helsinki, Department of Mathematics and Statistics, P.O. Box 68, 00014 University of Helsinki, Finland}
\email{jan.cristina@helsinki.fi}

\author[Iwaniec]{Tadeusz Iwaniec}
\address{Department of Mathematics, Syracuse University, Syracuse,
NY 13244, USA and Department of Mathematics and Statistics,
University of Helsinki, Finland}
\email{tiwaniec@syr.edu}

\author[Kovalev]{\\Leonid V. Kovalev}
\address{Department of Mathematics, Syracuse University, Syracuse,
NY 13244, USA}
\email{lvkovale@syr.edu}

\author[Onninen]{Jani Onninen}
\address{Department of Mathematics, Syracuse University, Syracuse,
NY 13244, USA}
\email{jkonnine@syr.edu}
\thanks{Cristina was supported by the Academy of Finland project 1123633 and ESF Network HCAA. Iwaniec was supported by the NSF grant DMS-0800416 and the Academy of Finland project 1128331. Kovalev was supported by the NSF grant DMS-0968756.
Onninen was supported by the NSF grant DMS-1001620.}




\begin{abstract} The central theme in this paper is
the Hopf-Laplace equation, which represents stationary solutions with respect to the inner variation of the Dirichlet integral.
Among such solutions are harmonic maps. Nevertheless, minimization of the Dirichlet energy among homeomorphisms often leads to mappings which are neither harmonic nor homeomorphisms. We prove that such mappings are harmonic outside of a singular set with small image. On the singular set they are locally Lipschitz, but not necessarily differentiable. 
\end{abstract}
\maketitle

 \tableofcontents

\section{Introduction}
 The  Hopf-Laplace equation arises in the study of the Dirichlet energy integral
\begin{equation}\label{direnergy}
\E_\X [h]= \iint_{\X} \abs{Dh}^2 =2 \iint_{\X} \left( \abs{h_z}^2+ \abs{h_{\bar z}}^2 \right)\, \dtext z
\end{equation}
for homeomorphisms $h \colon \X \to \Y$ between two designated domains $\X$ and $\Y$ in the complex plane $\C = \{z=x_1+ix_2 \colon x_1, x_2 \in \R\}$. In the recent paper~\cite{IKKO} we established some sufficient, and necessary conditions for the minimum of energy
to be attained by a harmonic diffeomorphism. In general, minimization of~$\E_\X$ yields non-harmonic, non-diffeomorphism, non-smooth solution of the Hopf-Laplace equation. The goal of this paper is to establish the regularity theory for such solutions. 

Throughout this text we  take advantage of the complex partial derivatives
\[
h_z  = \frac{\partial h}{\partial z}= \frac{1}{2} \left( \frac{\partial h}{\partial x_1} - i \frac{\partial h}{\partial x_2}  \right)  \quad \mbox{and} \quad
h_{\bar z}  = \frac{\partial h}{\partial \bar z}= \frac{1}{2} \left( \frac{\partial h}{\partial x_1} + i \frac{\partial h}{\partial x_2}  \right).
\]
In fact complex notation will be indispensable for advancing this work, especially when quadratic differentials will enter the stage.
Let us commence with the  variational formulation of the classical Dirichlet problem. One asks for the energy-minimal mapping $h \colon \X \to \C\,$ of the Sobolev class $h_\circ + \W_\circ^{1,2} (\X\!\shortrightarrow\!\C)\,$ whose boundary values are prescribed by means of a given mapping $h_\circ \in \W^{1,2} (\X\!\shortrightarrow\!\C)$. The \textit{first variation}, $\,h \leadsto h\,+ \,\epsilon \,\eta $,  in which $\eta \in \mathscr C^\infty_\circ (\X\!\shortrightarrow\!\C)$ can be any test function and $\epsilon \to 0$, leads to the Laplace equation
\begin{equation}\label{laplaceeq}
\Delta h = 4 h_{z \bar z} = 4 \frac{\partial^2 h}{\partial z \partial \bar z}=0.
\end{equation}

However, this approach is invalid when one seeks to minimize $\E_\X [h]$ among homeomorphisms; the injectivity of $h\,+\,\epsilon \,\eta$ is usually lost. The difficulty is circumvented by performing the  {\it inner variation}
\begin{equation}\label{innervar}
h \leadsto h \circ \chi_\epsilon\,, \qquad \chi_\epsilon \colon \X \onto \X
\end{equation}
where $\chi_\epsilon$ are $\mathscr C^\infty$-smooth automorphisms of $\X$ onto itself, defined for parameters $\epsilon \approx 0$, and $\chi_\circ= \mathrm{id} \colon \X \onto \X$. This is simply a change of the independent variable in the domain of definition of $h$. Here each $\chi_\epsilon$ may be the identity on $\partial \X$, but it need not be. The latter situation is called slipping along the boundary. The first derivative test $\frac{\dtext}{\dtext \epsilon} \E_\X [h \circ \chi_\epsilon]=0$, at $\epsilon =0$, yields
what we call the {\it Hopf-Laplace equation}
\begin{equation}\label{hopfeq}
\frac{\partial}{\partial \bar z} \big( h_z \overline{h_{\bar z}}\big)=0, \quad \mbox{in the sense of distributions.}
\end{equation}

The equation~\eqref{hopfeq} appears in the literature under several names. See, e.g., the analysis of ``variation a)'' in~\cite{Cob}, Chapter 1 of~\cite{Job} and the (more general)  energy-momentum equations in~\cite{SS, Ta}.

We shall investigate the Hopf-Laplace equation independently of its roots. Nevertheless,  because  of its affiliation to the Dirichlet integral the natural domain of definition, that we shall always assume here, sometimes implicitly,  is the Sobolev space $\W^{1,2}(\X)$. Thus the differential expression $h_z \overline{h_{\bar z}}$ will represent an integrable function. By Weyl's Lemma the equation reduces equivalently to the first order nonlinear PDE,
\begin{equation}\label{hopfeq2}
h_z \overline{h_{\bar z}}= \varphi\,, \quad \mbox{where } \varphi \mbox{ is a holomorphic function in $\X$.}
\end{equation}
We view $\varphi$ also as unknown quantity.
In this formulation one can speak of generalized solutions in the Sobolev space $\W^{1,1}_{\textnormal{loc}}(\X)$. But we shall study, predominantly, the \textit{orientation preserving mappings}, meaning that the Jacobian determinant is nonnegative almost
everywhere
\[
 J_h(z) = J(z,h) = |h_z|^2 - |h_{\bar{z}}|^2\;\geqslant 0.
\]
Within the class of orientation preserving mappings we do not capture generalized solutions. Indeed, in this case
{$ \abs{h_{\bar{z}}}^2 {=  \abs{h_{\bar{z}}}  \abs{h_{\bar{z}}} \le  \abs{h_{{z}}}  \abs{h_{\bar{z}}}  =} \abs{ h_z\,\overline{h_{\bar{z}}} } = \abs{\varphi }  \in \mathscr L^\infty _{\textnormal{loc}} (\mathbb X)$},
thus $h$ is (locally) a \emph{qc-deformation} in the sense of Ahlfors~\cite{Ah}. Since the complex derivatives $h_z$ and $h_{\bar z}$ are intertwined by the
Beurling-Ahlfors transform, it follows that $h_z  \in  BMO_{\rm loc}(\mathbb X)$.
In particular, for all $1\le  p  < \infty $ we have  $h \in \mathscr W^{1,p}_{\textnormal{loc}} (\mathbb X)$, hence
$h \in \mathscr C^\alpha_{\textnormal{loc}} (\mathbb X)$ for all $0 \leqslant \alpha < 1$.
The Lipschitz regularity is much harder to handle; this is the intricate part of our paper. Note that the qc-deformations of Ahlfors are not necessarily Lipschitz.

Returning to the variational approach to the Hopf-Laplace equation, if one admits slipping along the boundary, then the inner variations of the Dirichlet energy  should include diffeomorphisms $\chi_\epsilon \colon \X \onto \X$ that are  free on $\partial \X$.  Consequently, supplementary equations on $\partial\mathbb X$ emerge, which are best stated in terms of the {\it holomorphic quadratic differential}.
\begin{equation}\label{quaddiff}
\varphi (z) \, \dtext z^2 = h_z \overline{h_{\bar z}} \; \dtext z^2.
\end{equation}
It is called the  {\it Hopf differential} of $h$ in recognition of the related work of H.~Hopf~\cite{HopfH}.
Stated informally, the additional boundary equations~\cite[Lemma~1.2.5]{Job} say that $\varphi (z)\, \dtext z^2$ is real along $\partial \X$, see Definition~\ref{realqd}. This additional boundary equation will be satisfied by {\it minimal deformations}, see~\eqref{min}.
Clearly, conformal automorphisms $\chi_\epsilon \colon \X \onto \X$ do not change the energy. More generally, any conformal transformation of the domain $\mathbb X$ does not affect the Hopf-Laplace equation (\ref{hopfeq}).
In fact, because of that, it is the shape of the target $\mathbb Y$ and its closure, called \textit{deformed configuration}, that will really matter in questions to follow.

Naturally, complex harmonic functions solve the Hopf-Laplace equation. Worth noting, is that real-valued solutions in  $\mathscr C_{\loc}^1 (\X)$ must be harmonic.
If $\,h\,$  is $\,\mathscr C^2$-smooth then the Hopf-Laplace equation yields
\begin{equation}\label{HopfLapl}
 J(z,h) \, \Delta h = 0\;,\;\;\textnormal{where}\;\;  \Delta = 4 \frac{\partial^2}{\partial z \partial \bar {z}} \;\textnormal{ is the complex Laplacian.}
\end{equation}
Thus $\mathscr C^2$-solutions are harmonic in the region where the Jacobian determinant $J(z,h) \neq 0$.  There are other situations where the Hopf-Laplace equation implies harmonicity, see Theorem~\ref{thharm}. Nevertheless,   nonharmonic solutions arise naturally in global analysis of minimal surfaces~\cite[Ch. 2]{DHTb}, \cite{Ch1} and in the calculus of variations~\cite[Ch. 21]{AIMb}, \cite{Jono}. In mechanics, in particular in elasticity theory, the Hopf-Laplace equation appears under the name of  the energy-momentum equations.

Let us look at some  elementary though critical  examples. One might expect
from~\eqref{HopfLapl} that all $\mathscr W^{1,2}$-solutions whose Jacobian determinant does not vanish (almost everywhere) are harmonic maps.  This, however, is easily seen to be false, for in the complex plane the piece-wise linear mapping
\[
h(z) =
\begin{cases} 2z + \bar{z}\;, \quad &\textnormal{if} \;\;\;        \im  z\, \geqslant \,0 \\
 z + 2\bar{z} \;, \quad &\textnormal{if} \;\;\; \im z\, \leqslant \,0;
\end{cases}
\quad \quad \varphi = h_z\, \overline{h_{\bar{z}}}  \equiv 2\,,\quad\quad J(z, h) \equiv \pm 3
\]
is not harmonic. Observe that the Jacobian determinant changes sign.
More elaborate examples are provided in Section~\ref{examples:section},  which reveal that:

\begin{example}\label{exnoimp}  For each  $\,1< p < \infty\,$ there exists a generalized solution $\,h : \mathbb D \rightarrow \mathbb C\,$ (in the unit disk $\,\mathbb D\,$)  to the Hopf-Laplace equation
$h_z\, \overline{h_{\bar{z}}} \equiv 1$
 whose first derivatives belong to the Marcinkiewicz weak space $\, \mathscr L^p_{\textnormal {weak} } (\mathbb D)\,,$ but not to $\, \mathscr L^p (\mathbb D)$.
\end{example}

 The message from this example is that without supplementary conditions of a topological nature the Hopf-Laplace equation will not guarantee any substantial improvement of the regularity of the solutions. This is in marked contrast to the case of elliptic PDEs, like the  Laplacian.
 Consequently, we focus on solutions in a suitable closure of homeomorphisms.  Mathematical models of nonlinear elasticity motivate our calling such solutions  {\it {Hopf deformations}}.

The obvious question to ask is what  topological preconditions do we really need?  Partial answers are given in a few already known regularity results. H\'elein~\cite{He}  proved that quasiconformal solutions are harmonic. This is actually true~\cite{IKOhopf} for any homeomorphic solution in the Sobolev space $\,\mathscr W^{1,2}(\X\!\shortrightarrow\!\Y)\,$. Even more, injectivity is not necessary as long as the solution represents an open discrete map, such are the mappings of integrable distortion~\cite{IS}.
In particular, $\,\mathscr W^{1,2}$-solutions (not necessarily injective) with Jacobian $\, J(z,h) \geqslant \mbox{const} > 0 \,$ are also harmonic. However, in Example~\ref{exsqueeze}, we give a Lipschitz solution  with $\, J(z,h)  > 0 \,$, almost everywhere in the unit disk $\,\mathbb D\,$, which is not $\,\mathscr C^1$-smooth.  Even more, this solution is a homeomorphism in $\,\mathbb D\,,$ except for a tiny crack along the radius which is squeezed   into a point. Such a failure of injectivity along the cracks (reminiscent of elastic deformations) typically occurs when we pass to a weak limit of a minimizing sequence of homeomorphisms. With the loss of injectivity the extremal mappings will not be harmonic.  The examples, like the one mentioned above,  demonstrate  that the solutions can be Lipschitz at best. Our main results can be described as follows
\begin{itemize}
\item Partial harmonicity (Theorem~\ref{thharm}): Any Hopf deformation $h$ restricts to a harmonic diffeomorphism of $h^{-1}(\Y)$ onto $\Y$.
\item Lipschitz regularity (Theorem~\ref{MainTh}): Any Hopf deformation is locally Lipschitz.
\item Improved regularity (Theorem~\ref{smooth}): If $\Y$ is a polygon-type domain, then any $\mathscr C^1$-smooth minimal deformation is a harmonic diffeomorphism.
\item Non-smoothness of minimizers (Corollary~\ref{cor45}): minimizing deformations are not $\mathscr C^1$ in general.
\end{itemize}

In particular, our regularity results  apply to the minimal deformations; that is, to minimizers of the Dirichlet energy in the class of suitably generalized homeomorphisms.  An explicit example of a nonharmonic  minimal deformations is the following~\cite{AIM}:

\begin{example}\label{exhammering}  Consider two annuli $\mathbb X =  \{ z \colon \;r < |z| < R\}$ and $\mathbb Y =  \{ w \colon \, 1 < |w| < R_*\}$ where  $0<r<1<R < \infty$ and $R_*  = \frac{1}{2} ( R + R^{-1})$. The map
 \begin{equation}\label{hammering}
 h(z) =
\begin{cases}
  \;\frac{z}{|z|} \;, \;&\textnormal{if} \;\;\; r< |z| \leqslant 1 \;,\;\;\textnormal{squeezing to the unit cirle}
  \phantom{\bigg|} \\
   \;\frac{1}{2}(z + \frac{1}{\bar{z}} )  \;, \;&\textnormal{if} \;\;\; 1\leqslant |z| < R \;,\;\; \textnormal {the Nitsche harmonic map}
\end{cases}
\end{equation}
 takes $\,\mathbb X\,$ into $\,\overline{\mathbb Y}\,$. It satisfies the Hopf-Laplace equation
 \begin{equation}\label{Example2}
    h_z\, \overline{h_{\bar{z}}}  =  \varphi(z)  = \frac{- 1}{4 z^2}.
 \end{equation}
 The quadratic differential $\,\varphi(z)\, \textnormal{d} z^2\,$ is real and positive along the boundary circles.
\end{example}

\begin{center}
\begin{center}\begin{figure}[h]
\includegraphics[width=0.7\textwidth]{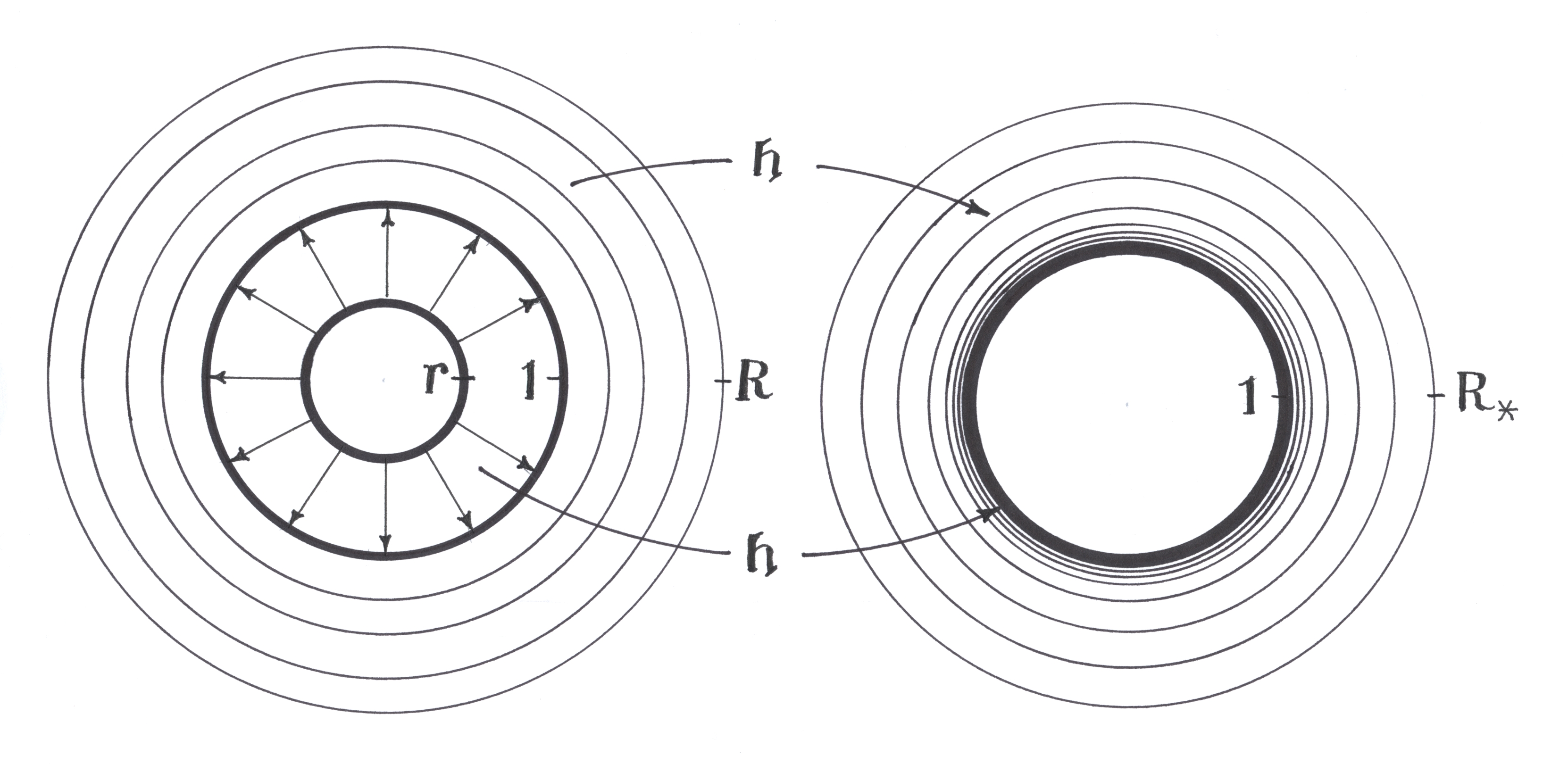}\caption{A nonharmonic minimal deformation}
\end{figure}\end{center}
\end{center}

The regularity theory of energy minimizing mappings has a long history. The monographs~\cite{Heb, LWb} provide an overview of the subjects. Note that in~\cite{Heb}  the solutions of the Hopf-Laplace equation are
called \emph{weakly Noether harmonic maps}. Numerous recent studies concern the Lipschitz continuity of energy minimizers in the setting of metric spaces, where higher degrees of smoothness are not
available~\cite{Ch,DM,GS,Jo96,Jo,KS,Ku,Wa}. However, in these studies minimization is performed among all mappings that  are either
prescribed on the boundary or belong to a given homotopy class. This minimization problem allows one to use the \textit{first variation}. Our approach is different
in that we obtain the Lipschitz continuity using only the inner variation. Thus we have in our disposal only the Hopf-Laplace equation instead of the Laplacian.

Before rigorous statements we need to review some notation and basic definitions.

\subsubsection{Domains}\label{sec101} We shall be concerned with mappings $\,h : \X \rightarrow \Y\,$ between bounded planar domains of finite connectivity $\, 1 \leqslant \ell < \infty\,.$ Thus each boundary $\,\mathfrak X = \partial \X\, $ and $\,\Upsilon = \partial \Y\, $ consists of $\,\ell\,$  disjoint continua.
 We reserve the notation,
\begin{equation}\label{boco}
\begin{split}
\mathfrak X_1, \mathfrak X_2, ..., \mathfrak X_\ell \;,\;\;\;\;\;\; &\textnormal{for the components of}  \;\;\mathfrak X= \partial \X  \\
\Upsilon_1, \Upsilon_2, ..., \Upsilon _\ell  \;,\;\;\;\;\; &\textnormal{for the components of} \;\; \Upsilon =\partial \Y.
\end{split}
\end{equation}

\subsubsection{Boundary Correspondence}

Every homeomorphism $\, h\colon \X  \onto \Y $ gives rise to a one-to-one correspondence between boundary components of $\X$ and boundary components of $\Y$.
It will involve no loss of generality in assuming (by re-arranging the indices, if necessary) that the correspondence is:
\begin{equation}\label{Corr}
h \colon \mathfrak X_\nu \rightsquigarrow \Upsilon_\nu \;,\;\;\; \textnormal{for}\;\;\; \nu = 1, 2, ..., \ell.
\end{equation}
This simply means that $h(x) \to \Upsilon_\nu$ as $x\to \mathfrak X_\nu$.
The above definition of the boundary correspondence is also pertinent to more general maps $\, h: \X \onto \Y$, specifically to those which are proper and
monotone in the sense defined below.

\begin{definition}~\cite{Mc, Rab, Wh}
A continuous mapping $f\colon \X \to \Y$ between metric spaces $\X$ and $\Y$  is {\it monotone} if for each $y\in f(\X)$ the set $f^{-1}(y)$ is compact and connected.
It is \emph{proper} if for each compact set $\mathbb F\subset f(\X)$ the set $f^{-1}(\mathbb F)$ is also compact.
\end{definition}

\subsubsection{$cd-$uniform convergence and the class $\mathscr H_{cd}(\X, \Y) $}
The idea of the $cd$-limit is a useful compromise between the concepts of $c$-uniform and uniform convergence.
 \begin{definition} A sequence of mappings $h_k: \X \rightarrow \Y \,,\, k = 1,2,...\,, $ is said to converge \textit{$cd$-uniformly}  to a mapping  $h: \X \rightarrow \mathbb R^2  $  if
 \begin{itemize}
 \item $h_k\rightarrow h$ \;\; $c$-uniformly \,(uniformly on compact subsets of $\X$)
 \item $\dist( h_k(x),\, \partial \Y) \; \rightarrow \;\textnormal{dist}( h(x),\, \partial \Y) $ \;\;uniformly  in $ \X\,.$
 \end{itemize}
We shall write it as  $h_k  \cto h\,$ and denote the class of $cd$-limits of homeomorphisms $\, h_{k}\colon \X  \onto \Y \,$
satisfying~\eqref{Corr} by $\mathscr H_{cd}(\X, \Y)$.
\end{definition}

 We emphasize that the boundary points of $\Y$ may be in the range of  $h \in\mathscr H_{cd}(\X, \Y)$. This fact
is crucial for several results that follow.  Precisely, we have
\begin{equation}\label{correspondence}
 \Y \subset h(\X) \subset \overline{\Y}  \;, \;\;\textnormal{for every }\;\;\; h \in\mathscr H_{cd}(\X, \Y).
 \end{equation}


If the range of a mapping $\;h \in\mathscr H_{cd}(\X, \Y)\;$ equals $\Y\,$ then $\,h\,$ is both  monotone and proper.
Actually, we have an even more precise statement.

\begin{proposition}[The Youngs approximation]\label{uniappthm} A continuous mapping
  $ h \colon \X  \to \Y $ between bounded $\ell$-connected  domains
belongs to $\mathscr H_{cd}(\X, \Y) $  if and only if it is monotone{,} proper and surjective.
  \end{proposition}

Let us introduce the class $\MPS(\X\!\shortrightarrow\!\Y)$ of continuous mappings $\,h \colon \X \to \Y\,$  between bounded domains, which are monotone,  proper and surjective. Note that such mappings,   as opposed to $\mathscr H_{cd}(\X, \Y) $, do not take values in $\,\partial\Y\,$.
Our notation is meant to emphasize this distinction.
The following is the extension of our earlier result~\cite{IKOhopf} on approximation  of homeomorphisms to the setting of monotone proper mappings.

\begin{theorem}[Approximation of $\mathcal{MPS}$ maps]\label{thm:diffeoapprox}
Let $h\colon \X\onto\Y$ be a continuous,  monotone {and} proper mapping of Sobolev class $\W^{1,2}_{\loc}(\X)$
between bounded $\ell$-connected domains $\X$ and $\Y$.  Then there exist diffeomorphisms
$h_{k}\colon \X\onto \Y$ such that
\begin{itemize}
\item $h_{k}- h\in  \mathscr A_\circ (\X)$
\item $\lim\limits_{k\to \infty}\|h_{k}-h\|_{\mathscr A(\X)}\; \;=0$.
\end{itemize}
\end{theorem}

Hereafter $\mathscr A (\X)= \mathscr C (\overline{\X}) \cap \W^{1,2} (\X\!\shortrightarrow\!\C)$ is the Royden algebra equipped with the norm
\begin{equation}\label{Royden algebra}
 \| f\|_{\mathscr A(\X)} = \| f\|_{\mathscr C(\X)} \;+\; \| Df\|_{\mathscr L^2(\X)}.
\end{equation}
and $\mathscr A_\circ (\X)  = \mathscr C_\circ(\overline{\X}) \cap \W_\circ^{1,2}(\X\!\shortrightarrow\!\C)\,$ is a subalgebra obtained by completing the space  $\mathscr C_{\circ}^\infty (\X)$ with respect to this norm.

The following  corollary is immediate from Theorem~\ref{thm:diffeoapprox}.

\begin{corollary}\label{jacdncs}
If $h \in \MPS(\X\!\shortrightarrow\!\Y) \cap \W^{1,2}_{\loc}(\X)$, then the Jacobian of $h$ does not change sign; that is, either $J_h \ge 0$ a.e. or $J_h \le 0$ a.e.
\end{corollary}

\subsubsection{Deformation and the class $\,\mathfrak D(\X\,,\Y)\,$}
The concept of \textit{deformation} seems to differ only a little from the routinely used $\,\mathscr W^{1,2}\,$--weak limit of  homeomorphisms.
However, when the topological features of the minimal-energy solutions are of major concern the \textit{deformations} become better suited than the weak limits of homeomorphisms.
For example, in~\cite{IKKO} the concept of  \textit{deformations} was critical in establishing existence of harmonic homeomorphisms between planar doubly connected
domains.  Before making the precise definition let us look at some examples  of both geometric and analytical nature.

 A sequence of homeomorphisms converging weakly in $\,\mathscr W^{1,2}\,$  actually converges $c$-uniformly, so the limit is a continuous map. But that is  all, in general,  what the weak limit of a minimizing sequence can receive from homeomorphisms. Consider the classical example:
 \begin{example}
 Conformal automorphisms $\, h_k \colon \mathbb D  \onto \mathbb D\, $ of the unit disk, $\, h_k(z) = \frac{z- a_k}{1 - z \overline{a_k}}\,,\; |a_k| \rightarrow 1 \,,$ are the minimal energy maps. They  converge  weakly in $\,\mathscr W^{1,2}(\mathbb D)\,$ and $\,c$-uniformly  to a constant map.  Such a trivial loss of topological distinctions is due to the failure of  $\,cd$-convergence.
 \end{example}
 On the other hand, when Sobolev mappings  come into play, we find that a $c$-uniform limit of homeomorphisms  $\, h_k \colon \X  \onto \Y\,$ that are  uniformly  bounded in $\;\mathscr W^{1,2}(\X\!\shortrightarrow\!\Y)\,$   satisfies the measure theoretical condition of \textit{non-overlapping}
 \begin{equation}\label{nonover}
  \iint_\X | J(z,h)| \,\textnormal{d}z  \leqslant |\Y| ,
  \end{equation}
  due to the $\mathscr L^1$-weak subconvergence of Jacobians~\cite[Theorem 8.4.2]{IMb}.
 In nonlinear elasticity theory this may be interpreted as saying that \textit{interpenetration of matter} does not occur. However, we will be forced to take into consideration  more general minimizing sequences
 $h_k  \overset{\textnormal{\tiny{cd}}}{\longrightarrow} h\,$ in which the mappings are neither homeomorphisms nor they have uniformly bounded energy.
We still impose the nonoverlapping condition ~\eqref{nonover}.
In view of Corollary~\ref{jacdncs} the Jacobian of $h$ does not change sign, and thus it entails no loss of generality to assume that $J_h\ge 0$ a.e.
These observations drive us to the following

\begin{definition}[Deformation] \label{deformations} Let $\X$ and  $\Y$ be bounded $\ell$-connected domains. A deformation  is a mapping  $h \in \mathscr H_{cd}(\X, \Y)$ such that
 \begin{itemize}
 \item  $\, h \in \mathscr W^{1,2} (\X\!\shortrightarrow\!\mathbb R^2)$
 \item The Jacobian determinant of $h$ is nonnegative a.e. and satisfies the  \textit{non-overlapping condition }
 \begin{equation}\label{nonoverlapping}
 \iint_\X  J(z,h) \,\textnormal{d}z  \leqslant \abs{\Y}.
 \end{equation}
 \end{itemize}
\end{definition}

The class of all deformations  will be denoted by $\mathfrak D(\X, \Y)$. Throughout what follows, if no confusion can arise,  we shall freely assume
without explicit mention that the class  $\mathfrak D(\X,\Y)$ is nonempty. We again strongly emphasize that a deformation $h \in \mathfrak D(\X,\Y)$ may take points
of $\X$ into $\partial\Y$, a key point that will affect the forthcoming arguments.  The structure of the preimage $h^{-1}( y_\circ)$ of a point
in $y_\circ \in \partial \Y$ is a delicate issue that we address in sections~\ref{secpreimages} and~\ref{secsmo}.

In~\cite{IKKO} we have already established the essential properties of deformations. In particular, the class $\mathfrak D(\X\,,\Y)$ is sequentially weakly closed in  $\mathscr W^{1,2} (\X\!\shortrightarrow\!\mathbb R^2)$, provided that $\X$ has $2 \le \ell < \infty$ boundary components none of which are points~\cite[Lemma 3.13]{IKKO}. Under these hypotheses we have the following.
\begin{corollary} The infimum energy of the Dirichlet energy within the class  $\mathfrak D(\X\,,\Y)$ is attained.
\end{corollary}

A mapping $h_\circ \in  \mathfrak D (\X\,,\Y)  $ such that
\begin{equation}\label{min}
\int_\X \abs{Dh_\circ}^2 \, \textnormal dx = \min_{h\in\mathfrak D(\X\,,\Y)} \int_\X \,|\,Dh(x)\,|^2 \; \textnormal d x
\end{equation}
will be called a \textit{minimal deformation}.

 \subsection{Hopf deformations}

\begin{definition}\label{hopfdefdef}
 Let $\X$ and $\Y$ be bounded $\ell$-connected domains.
 A deformation $\, h \in \mathfrak D (\X\,,\Y)$ which satisfies the Hopf-Laplace equation
\begin{equation}\label{hlagain}
\frac{\partial}{\partial \overline{z}} \left(\, h_z \, \overline{h_{\overline{z}}} \,\right)  =  0
\end{equation}
is called a {\it Hopf deformation}.
\end{definition}
Thus Hopf deformations are among stationary solutions (with respect to the inner variation) of the Dirichlet energy integral. We shall prove that,

\begin{theorem}\label{thharm} Any Hopf deformation $h\in\mathfrak D(\X,\Y)$ is  a harmonic diffeomorphism of $h^{-1}(\Y) \subset \X$ onto $\Y$.
\end{theorem}

Observe that in Example~\ref{exhammering} the mapping $h$ fails to be harmonic exactly in the part of $\X$ that is squeezed into a boundary component of the target annulus $\Y$. More precisely, this particular component of $\partial\mathbb Y$ is the inner boundary circle at which $\mathbb Y$ is not convex.
In particular, when the range of $h$ is $\mathbb Y$, we can combine Theorem~\ref{thharm} with Proposition~\ref{uniappthm} to obtain

\begin{corollary} Let $\X$ and $\Y$ be bounded $\,\ell$-connected domains and
 $ h\colon \X  \overset{\textnormal{\tiny{onto}}}{\longrightarrow} \Y$ a monotone {and} proper mapping of Sobolev class $\mathscr W^{1,2}(\X\!\shortrightarrow\!\Y)$ that satisfies the
 Hopf-Laplace equation~\eqref{hlagain}.
Then $h$ is a harmonic diffeomorphism.
\end{corollary}

This corollary represents another advance on the Eells-Lemaire problem~\cite{EL2, EL3}: under what conditions does the holomorphicity of the Hopf differential
imply that the mapping is harmonic?

\subsection{Lipschitz regularity}

The foremost interesting and much harder task is to determine the regularity of a Hopf deformation when its image $h(\X)$ goes over the designated target $\Y$; that is,   into its closure. There are minimal deformations, in particular Hopf deformations, which are not $\mathscr C^1$-smooth in $\X$, see Corollary~\ref{cor45}. Our main result in this paper is the following best possible regularity of Hopf deformations.

\begin{theorem} \label{MainTh} 
Every Hopf deformation $h \in \mathfrak D(\X\,,\Y)$
is locally Lipschitz continuous.  In particular, a minimal deformation of planar $\,\ell$-connected domains   is locally Lipschitz continuous.
\end{theorem}

\begin{theorem}\label{smooth} 
Suppose that $ h\in \mathfrak D(\X,\Y) $ is a minimal deformation and the nonconvex part of $\partial \Y$ is at most countable.
Then $h$ is $\mathscr C^1$-smooth if and only if it is a harmonic homeomorphism.
\end{theorem}

The term \textit{nonconvex part of} $\partial \Y$ refers to a set of points in $\partial \Y$ where $\Y$  is not  convex in any neighborhood. Precisely,
\begin{definition}
We say that $\Y$ is convex at $y_\circ \in \partial \Y$ if for some $\varepsilon > 0$ the set $\{ y \in \Y \,: \, |y - \,y_\circ | < \varepsilon\,\}$ is convex.
\end{definition}
A natural example of a domain with finite nonconvex part of the boundary is an $\ell$-connected polygonal domain.

According to~\cite[Theorem 2.4]{IKKO} there exists a nondecreasing function
\[
\Theta \colon (0,\infty)\to (0,\infty), \qquad \lim\limits_{\tau\to\infty}\frac{\Theta(\tau)}{\tau}= 1
\]
such that the following holds. Whenever two bounded doubly connected domains $\X$ and $\Y$ admit an energy minimizing diffeomorphism $h\colon \X \onto\Y$,
we have
\begin{equation}\label{eme1}
\Mod\Y \;\ge\; \Theta(\Mod\X)
\end{equation}
where $\Mod$ stands for the conformal modulus. Now combining Theorem~\ref{smooth} with the nonexistence of energy minimizing diffeomorphisms, we arrive at the following corollary.

\begin{corollary}\label{cor45}
Let $\X$ and $\Y$ be bounded doubly-connected domains such that the nonconvex part of $\partial \Y$ is at most countable
and~\eqref{eme1} fails. Then any minimal deformation in  $\mathfrak D(\X,\Y)$ is not $\mathscr C^1$-smooth.
\end{corollary}

\begin{center}\begin{figure}[h]
\includegraphics[width=0.8\textwidth]{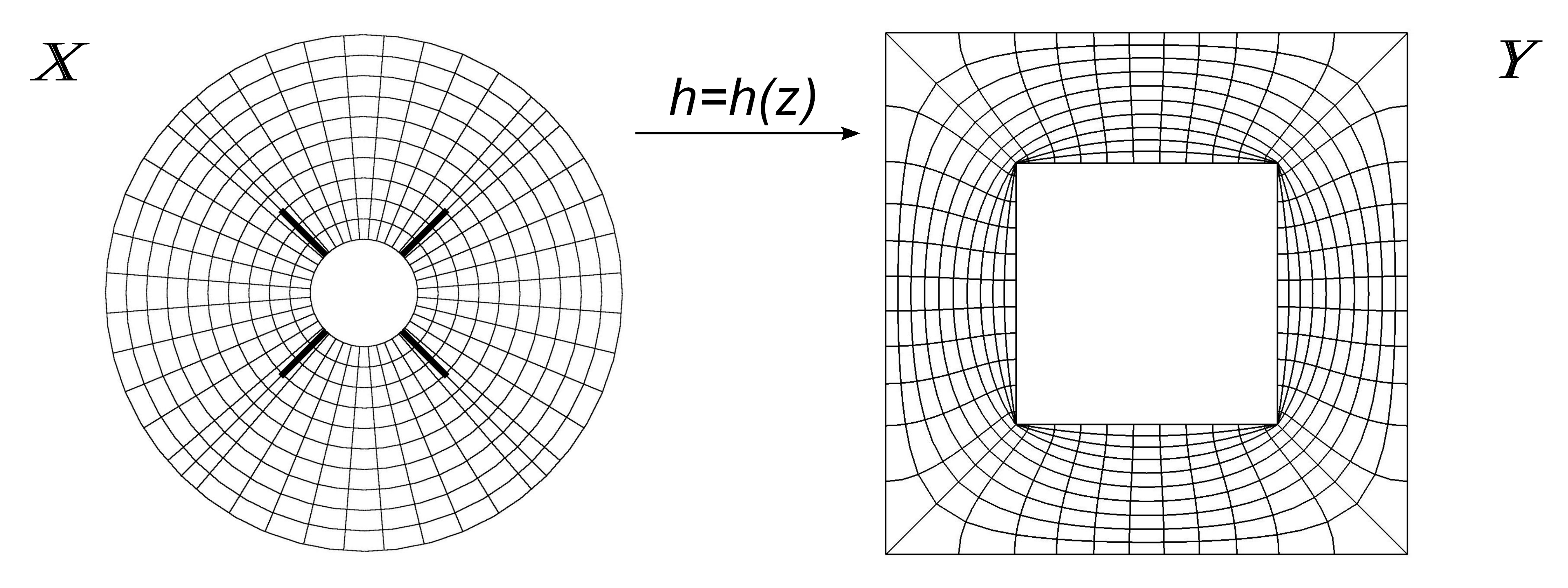} \caption{A minimal deformation fails to be $\mathscr C^1$-smooth.}\label{fig4}
\end{figure}\end{center}

We conclude this introduction with a brief outline of the proofs of the main results. To prove the partial harmonicity of a Hopf deformation $h$ we show that
$h$ has a stronger energy minimization property in $h^{-1}(\Y)$ than in all of $\X$, namely~\eqref{goal} holds. The three main
ingredients of this proof are: diffeomorphic approximation of monotone Sobolev mappings (\S\ref{secapp}),
a Reich-Strebel-type inequality (Lemma~\ref{reich-strebel}), and the structure of preimages of points under $h$ (\S\ref{secpreimages}).
The partial harmonicity of $h$ is essential in the proof of its Lipschitz continuity, however a different approach
is required for the set $\X\setminus h^{-1}(\Y)$. To this end we use a Bonnesen-type inequality, i.e., a stability result for the isoperimetric
inequality (\S\ref{sectionlip}). The proof of Theorem~\ref{smooth} invokes Besicovitch's removability theorem for holomorphic functions, further analysis of
preimages of points under $h$, and the boundary point lemma of E.~Hopf.

One method of proving Lipschitz continuity of $h$ proceeds through the subharmonicity of the energy density $\abs{Dh}^2$, as in
~\cite{GS} and~\cite{Ch}. This method relies on the first variation of $h$, which is not available to us.
Although the Hopf-Laplace equation appears to imply the subharmonicity of $\abs{Dh}^2$ on the formal level, our Example~\ref{exnoimp} shows that this is not the case.

\section{Preliminaries}\label{secpre}



The uniform limit of self-homeomorphisms of the sphere $\s^2$ is a monotone surjective mapping~\cite[IX.3.11]{Wh}. In the converse direction,
a monotone mapping of $\s^2$ onto itself can be refined in some Jordan subdomains of $\s^2$ in which it becomes a homeomorphism.

\subsection{The Youngs refinement}
Let $\X$ and $\Y$ be $\ell$-connected domains and $f\in \MPS(\X\!\shortrightarrow\!\Y)$.
Recall that $\mathbb D \Subset \Y$ is a Jordan domain if it is the interior of a homeomorphic image of the closed unit disk.

\begin{proposition}\label{homext}
For every Jordan domain $\mathbb D \Subset \Y$ and $f\in \MPS(\X\!\shortrightarrow\!\Y)$ the preimage $\mathbb U= f^{-1}(\mathbb D) \Subset \X$ is a simply connected domain in $\X$.
Furthermore, there is $f_{\mathbb D} \in \MPS(\X\!\shortrightarrow\!\Y)$, {\it the Youngs refinement}, such that
 \begin{itemize}
\item ${f}_{\mathbb D}=f$ in $\X \setminus \mathbb U$.
\item $f_{\mathbb D}$ restricted to $\mathbb U$ is a homeomorphism of $\mathbb U$ onto $\mathbb D$.
\end{itemize}
\end{proposition}
\begin{proof}
This is a consequence of the Youngs modification theorem~\cite[Theorem 10.1]{Yo} or ~\cite[II.1.47]{Rab} about continuous monotone mappings $F\colon \s^2 \onto \s^2$.
It says that to every Jordan domain $\Delta \subset \s^2$ there corresponds a continuous monotone mapping $F_{\Delta} \colon \s^2 \onto \s^2$ such that
\begin{itemize}
\item $F_{\Delta}=F$ in $\s^2 \setminus F^{-1} (\Delta)$
\item $F_{\Delta}$ maps $F^{-1} (\Delta)$ homeomorphically onto $\Delta$.
\end{itemize}
We prove Proposition~\ref{homext} by applying the Youngs theorem to the $\ell$-point compactifications of $\X$ and $\Y$.
Indeed, both $\X$ and $\Y$ are homeomorphic to the sphere $\s^2$ with $\ell$ punctures, say via the maps
$\Phi\colon \X\to \s^2\setminus\{x_1,\dots,x_{\ell}\}$ and  $\Psi\colon \Y\to \s^2\setminus\{y_1,\dots,y_{\ell}\}$
where the punctures are enumerated in the same way as the boundary components~\eqref{boco}.
For any mapping $f\in \MPS(\X\!\shortrightarrow\!\Y)$ the composition $\Psi \circ f \circ \Phi^{-1}$ extends to a continuous mapping
\begin{equation}\label{eq333}
F \colon \s^2 \onto \s^2, \quad \mbox{ with } F \big(\Phi (x_\nu)\big) = \Psi (y_\nu), \quad \nu=1, \dots , \ell.
\end{equation}
Any simply connected domain $\mathbb  D \subset \Y$ is mapped via $\Psi $ onto a simply connected domain $\Delta = \Psi (\mathbb D) \subset \s^2$
which stays away from the punctures. By the Youngs theorem $F^{-1} (\Delta)$ is a simply connected domain in $\s^2$.
It follows that $f^{-1} (\mathbb D)$ is a simply connected domain compactly contained in $\X$. The homeomorphism
$F_{\Delta} \colon F^{-1}(\Delta) \onto \Delta$ (the Youngs refinement) yields a homeomorphism $f_{\mathbb D} \colon f^{-1} (\mathbb D) \onto \mathbb D$.
\end{proof}

\subsection{The Youngs approximation, proof of Proposition~\ref{uniappthm}}

\begin{proof}
First, suppose that $f \colon \X \into \Y$ is a $cd$-limit of homeomorphisms $f_k \colon \X \onto \Y$. Therefore, $f(\X)=\Y$, by~\eqref{correspondence}.
The preimage of a compact set in $\Y$, under $f$, is closed and stays away from $\partial \X$
because $f$ satisfies~\eqref{Corr}. Thus $f$ is proper. To show monotonicity we appeal to the induced
mappings $F_k \colon \s^2 \onto \s^2$, $F_k = \Psi \circ f_k \circ \Phi^{-1}$ as in~\eqref{eq333}. These are homeomorphisms of $\s^2$ onto
itself converging uniformly to $F= \Psi \circ f \circ \Phi^{-1}$. Thus $F$ is monotone~\cite[Theorem 11.1]{Yo}.  Since $f \colon \X \onto \Y$ we
see that $f$ is monotone as well.

In the converse direction we must show that every mapping $f\in \MPS(\X\!\shortrightarrow\!\Y)$ belongs to $\Ho_{\cd} (\X, \Y)$. In the proof of Proposition~\ref{uniappthm} we shall appeal to the classical
Youngs approximation theorem~\cite[Theorem 11.1]{Yo} (or~\cite[II.1.57]{Rab}). It asserts that
\vskip0.2cm
{\it A continuous mapping $f \colon \s^2 \onto \s^2$  is monotone if and only if it is a uniform limit of homeomorphisms of $\s^2$ onto $\s^2$.}
\vskip0.2cm
Let us  compactify $\X$ and $\Y$ as in the proof of Proposition~\ref{homext}, via the homeomorphisms $\Phi \colon \X \onto \s^2\setminus\{x_1,\dots,x_{\ell}\}$ and
$\Psi \colon \Y \onto \s^2\setminus\{y_1,\dots,y_{\ell}\}$.
For any $f\in \MPS(\X\!\shortrightarrow\!\Y)$ the induced  mapping $F \colon \s^2 \onto \s^2$, being monotone, can be uniformly approximated by
homeomorphisms $F_k \colon \s^2 \onto \s^2$. By altering $F_k$ near punctures we can ensure that $F_k(x_\nu)=y_\nu$ for $\nu=1,\dots,\ell$.
We then return to the domains $\X$ and $\Y$ via $\Phi$ and $\Psi$ to define homeomorphisms
\begin{equation}\label{sphere}
f_k= \Psi^{-1} \circ F_k \circ \Phi \colon \X \onto \Y
\end{equation}
and observe that $f_k\cto f$ by construction.
\end{proof}

The Youngs refinement provides a homeomorphic replacement of a monotone mapping. The following proposition  shows that such a replacement can be chosen to be  harmonic. It combines classical results of potential theory~\cite{GMb, Ranb} with a recent extension of the Rad\'o-Kneser-Choquet theorem.

\begin{proposition}\label{poisson}
Let $\X\subset \C$ be a domain. To every bounded simply connected domain $\mathbb U \Subset \X$ there corresponds a unique linear operator
\[{\bf P}_{\mathbb U} \colon \mathscr C (\X) \to \mathscr C (\X)\]
such that for every $f\in \mathscr C(\X)$
\begin{enumerate}
\item[(i)] ${\bf P}_{\mathbb U} f =f$ in $\X \setminus \mathbb U$.
\item[(ii)] ${\bf P}_{\mathbb U} f$ is harmonic in $\mathbb U$.
\end{enumerate}
Such an operator has the additional properties ~\cite{IKKO}
\begin{enumerate}
\item[(iii)] ${\bf P}_{\mathbb U } f \in f+ \mathscr W^{1,2}_\circ (\mathbb U)$, whenever $f\in \mathscr C (\X) \cap \mathscr W^{1,2}_{\loc} (\X\!\shortrightarrow\!\C)$. Moreover,
\[\E_{\mathbb U} [{\bf P}_{\mathbb U} f] \le \E_{\mathbb U} [ f], \qquad \mbox{provided} \quad \E_{\mathbb U} [ f] < \infty .\]
\item[(iv)] Suppose that the restriction $f_{|_{\mathbb U}}$ of $f\in \mathscr C (\X)$ is a homeomorphism of $\mathbb U$ onto a convex domain $\mathbb D \subset \C$, then ${\bf P}_{\mathbb U} f \colon \mathbb U \onto \mathbb D$ is a harmonic diffeomorphism.
\end{enumerate}
\end{proposition}

We call ${\bf P}_\U f$ the \emph{harmonic replacement} of $f$.


\subsection{Prerequisites from holomorphic quadratic differentials}\label{secqd}
In this section we recall some basic facts about holomorphic quadratic differentials. The general reference for these topics is~\cite{Stb}.

Let $\X$ will be a bounded $\ell$-connected domain and $\varphi \colon \X \to \C$ a  holomorphic function with isolated zeros, called {\it critical points}.
Denote $\X_\circ=\X\setminus\{\text{zeros of }\varphi\}$.
In a neighborhood of every point $a\in \X_\circ$ one can introduce a local conformal mapping $w= \Phi (z)= \int \sqrt{\varphi (z)}\, \dtext z$, called a natural parameter near $a$. Through every regular point there pass two $\mathscr C^\infty$-smooth orthogonal arcs, called horizontal and vertical arcs.  A vertical arc is a $\mathscr C^\infty$-smooth curve $\gamma \colon t \to \gamma (t)$, $a<t<b$, along which
\[[\dot{\gamma} (t)]^2 \varphi \big(\gamma (t)\big)<0, \quad a<t<b.\]
A horizontal arc is a $\mathscr C^\infty$-smooth curve $\beta \colon t \to \beta (t)$, $c<t<d$, along which
\[[\dot{\beta} (t)]^2 \varphi \big(\beta (t)\big)>0, \quad c<t<b.\]
We emphasize that this yields, in particular, that such arcs only contain regular points of $\varphi$. A {\it vertical trajectory} of $\varphi$ in $\X$ is a maximal vertical arc; that is, not properly contained in any other vertical arc. Hereafter, with the customary abuse of notation, the same symbol $\gamma$ will be used for both the parametrization $\gamma = \gamma (t)$ and its range. Similarly, a \emph{horizontal trajectory} is a maximal horizontal arc in $\X$. Through every regular point of $\varphi$ there passes a unique vertical (horizontal) trajectory. A trajectory whose closure contains a critical point of $\varphi$ is called a {\it critical trajectory}. There are at most a countable   number of critical trajectories, so they cover a set in $\X$ of measure zero. We will be largely concerned with noncritical trajectories.

\begin{definition}($\varphi$-rectangle)
A $\varphi$-rectangle of a quadratic differential $\varphi (z)\, \dtext z^2$ is any simply connected domain $\mathfrak R \subset \X_\circ$ on which the natural parameter $w=\Phi (z)= \int \sqrt{\varphi (z)}\, \dtext z$ has a univalent branch which takes $\mathfrak R$ onto a Euclidean rectangle
\[\Phi (\mathfrak R)= \{w=t+i\tau \colon 0<t<T \quad \mbox{and} \quad a< \tau<b\}.\]
Note that $\mathfrak R$ contains no zeros of $\varphi$. We will be concerned with $\varphi$-rectangles which are compactly contained in $\X_\circ$, so $\Phi$ defines a diffeomorphism of a neighborhood of $\overline{\mathfrak R}$ onto a neighborhood of $\overline{\Phi (\mathfrak R)}$. Then we define the horizontal edges of $\mathfrak R$, $\alpha= \Phi^{-1} \left([0, T] \times \{a\}\right)$ and $\beta= \Phi^{-1} \left([0, T] \times \{b\}\right)$ and similar for the vertical edges.
\end{definition}

{\it Every noncritical vertical trajectory $\gamma \subset \Omega$  in a simply connected domain $\Omega$ is a cross cut}, see Theorem~15.1 in~\cite{Stb}. Thus  in the maximal interval $a<t<b$ of the existence of $\gamma = \gamma (t)$ both limit sets at the end-points of $\gamma$, denoted by $\gamma \{a\}$ and $\gamma \{b\}$, lie in $\partial \Omega$. Let $\gamma_\circ \subset \gamma$ be any closed vertical subarc of $\gamma$, defined by $\gamma_\circ (t)= \gamma (t)$ for $a_\circ \le t \le b_\circ$, where $a<a_\circ < b_\circ <b$. Then the $\varphi$-length of $\gamma_\circ$ is equal or smaller than $\varphi$-length of any rectifiable curve  $\beta \subset \Omega$ which connects $A_\circ = \gamma (a_\circ)$ with $B_\circ = \gamma (b_\circ)$. This means that
\begin{equation}\label{length0}
\int_{\gamma_\circ} \abs{\varphi}^{\nicefrac{1}{2}}\, \abs{\dtext z} \le \int_{\beta} \abs{\varphi}^{\nicefrac{1}{2}}\, \abs{\dtext z}
\end{equation}
see \cite[Theorem~16.1]{Stb}. Note that $\gamma \setminus \gamma_\circ$ consists of two components (two disjoint vertical arcs). Inequality~\eqref{length0} can be slightly generalized; it is not necessary to assume that the end-points of $\beta$ coincide with the endpoints of $\gamma_\circ$.

\begin{lemma}\label{lengthineq}
Let $\beta \subset \Omega$ be a locally rectifiable arc in a simply connected region whose closure intersects both components of $\gamma \setminus \gamma_\circ$, then
\begin{equation}\label{length}
\int_{\gamma_\circ} \abs{\varphi}^{\nicefrac{1}{2}}\, \abs{\dtext z} \le \int_{\beta} \abs{\varphi}^{\nicefrac{1}{2}}\, \abs{\dtext z}.
\end{equation}
\end{lemma}

\begin{proof}
Let $A, B \in \gamma$ be points in different components of $\gamma \setminus \gamma_\circ$ that are approachable through the arc $\beta$; that is,
\[A = \lim_{n \to \infty} A_n \quad \mbox{and} \quad B = \lim_{n \to \infty} B_n, \qquad \mbox{where } A_n, B_n \in \beta. \]
Let $[A,B]_\gamma$ denote subarc of $\gamma$ that connects $A$ and $B$. We certainly have $\gamma_\circ \subset [A, B]_\gamma$, so
\[  \int_{\gamma_\circ} \abs{\varphi}^{\nicefrac{1}{2}}\, \abs{\dtext z} \le \int_{[A,B]_\gamma} \abs{\varphi}^{\nicefrac{1}{2}}\, \abs{\dtext z}.
 \]
Similarly, we denote by $[B_n , A_n]_\beta \subset \beta$ the closed (rectifiable) subarc of $\beta$ which connects $B_n$ and $A_n$. Since $A_n \to A \in \gamma \subset \Omega$ and  $B_n \to B \in \gamma \subset \Omega$ for sufficiently large $n$ the straight segments $[A, A_n]$ and $[B_n, B]$ lie in $\Omega$. We now have a rectifiable curve  $[A, A_n] \cup [A_n , B_n]_\beta \cup [B_n, B]$ in $\Omega$ which connects the end-points of $[A,B]_\gamma$. Therefore,  we have
\[
\begin{split}
\int_{\gamma_\circ} \abs{\varphi}^{\nicefrac{1}{2}} & \le \int_{[A,B]_\gamma} \abs{\varphi}^{\nicefrac{1}{2}}  \le  \int_{[A_n,B_n]_\gamma} \abs{\varphi}^{\nicefrac{1}{2}} + \int_{[A,A_n]} \abs{\varphi}^{\nicefrac{1}{2}}   + \int_{[B_n,B]} \abs{\varphi}^{\nicefrac{1}{2}} \\
& \le  \int_{\beta} \abs{\varphi}^{\nicefrac{1}{2}}+ \left( \abs{A_n -A} + \abs{B_n-B} \right)\norm{\varphi}_{\mathscr C (\Omega)}^{\nicefrac{1}{2}} \longrightarrow \int_{\beta} \abs{\varphi}^{\nicefrac{1}{2}},
\end{split}
\]
as desired.
\end{proof}




The next lemma deals with a holomorphic quadratic differential $\varphi\, \dtext z^2$ which is real on the boundary of a $\mathscr C^1$-smooth domain $\X$, no single points as components of $\partial \X$.
\begin{definition}\label{realqd}
A quadratic differential $\varphi\, \dtext z^2$  is said to be real on the boundary of a $\mathscr C^1$-smooth $\ell$-connected domain $\X$ if  $\varphi$ is
smooth up to $\partial\X$ and each component of $\partial \X$ is either a horizontal or a vertical trajectory of $\varphi\, \dtext z^2$.
\end{definition}

\begin{lemma}\label{callsomething}
Let $\X$ be a finitely connected domain with $\mathscr C^1$-smooth boundary. Let
$\varphi\, \dtext z^2$ be a holomorphic quadratic differential in $\X$ which is real on $\partial \X$.
Suppose that $\Gamma$ is a vertical trajectory of $\varphi\, \dtext z^2$ with both ends approaching the same boundary component of $\X$. Then
the  components of $\X\setminus \Gamma$ are not simply connected.
\end{lemma}

\begin{proof} Suppose that the set $\X\setminus \Gamma$ has a simply connected component $\mathbb G$.
There are no closed trajectories in $\mathbb G$, for such a trajectory must enclose a pole of $\varphi$.
The global structure of trajectories of a holomorphic quadratic differential with finite norm~\cite{Stb} is inconsistent with $\partial \mathbb G$
being a union of a vertical trajectory and another (vertical or horizontal) trajectory.
\end{proof}

\begin{lemma}[Fubini-like integration formula]\label{Fublem}
Let $\varphi (z)\, \dtext z^2$ be a holomorphic quadratic differential in a simply connected domain $\Omega \subset \C$, $\varphi \not\equiv 0$. Suppose that $F$ and $G$ are measurable functions in $\Omega$ such that
\begin{equation}\label{eqfub1}
\iint_\Omega \abs{\varphi (z)} \abs{F(z)}\, \dtext x \dtext y < \infty \quad \mbox{and} \quad \iint_\Omega \abs{\varphi (z)} \abs{G(z)}\, \dtext x \dtext y < \infty.
\end{equation}
Then for almost every vertical trajectory $\gamma$ of $\varphi (z)\, \dtext z^2$ we have
\begin{equation}\label{eqfub2}
\int_\gamma \abs{\varphi (z)}^{\nicefrac{1}{2}} \abs{F(z)}\, \abs{\dtext z} < \infty \quad \mbox{and} \quad \int_\gamma \abs{\varphi (z)}^{\nicefrac{1}{2}} \abs{G(z)} \abs{\dtext z} < \infty.
\end{equation}
If, in addition,
\begin{equation}\label{eqfub3}
\int_\gamma \abs{\varphi (z)}^{\nicefrac{1}{2}} {F(z)}\, \abs{\dtext z} =\int_\gamma \abs{\varphi (z)}^{\nicefrac{1}{2}} {G(z)}\, \abs{\dtext z},
\end{equation}
then
\begin{equation}\label{eqfub4}
\iint_\Omega \abs{\varphi (z)} {F(z)}\, \dtext x \dtext y= \iint_\Omega \abs{\varphi (z)} {G(z)}\, \dtext x \dtext y.
\end{equation}
\end{lemma}

\begin{proof}
According to~\cite[\S19.2]{Stb} $\Omega$ can be covered, up to a set of measure zero, by a countable number of disjoint $\varphi$-strips. These are open connected subsets of $\Omega$ with no critical points, such that a locally defined analytic function $\Phi (z)= \int \sqrt{\varphi (z)}\, \dtext z$ is actually a univalent conformal mapping of the $\varphi$-strip onto a Euclidean vertical strip $\mathcal S$ in the $w$-plane, $w= \Phi (z)$.
\[\mathcal S = \{w=t+i\tau \colon 0<t<T, \quad \alpha (t)< \tau < \beta (t)  \}\]
where $-\infty \le \alpha (t) < \beta (t) \le + \infty$ are measurable functions. Making a substitution $z= \Phi^{-1}(w)$ the problem reduces equivalently to the usual Fubini's theorem in Euclidean vertical strip, for functions $\tilde{F}(w)= \abs{\varphi (z)} F(z)$ and  $\tilde{G}(w)= \abs{\varphi (z)} G(z)$.
\end{proof}

\begin{corollary}\label{fubcor}
Assume, instead of condition~\eqref{eqfub3} in Lemma~\ref{Fublem} that
\[\int_\gamma \abs{\varphi}^{\nicefrac{1}{2}} \abs{F} \le \int_\gamma \abs{\varphi}^{\nicefrac{1}{2}} \abs{G}
\]
for almost every noncritical vertical trajectory of $\varphi \, \dtext z^2$. Then
\[\iint_\Omega \abs{\varphi}\, \abs{F} \le \iint_\Omega \abs{\varphi}\, \abs{G} .\]
\end{corollary}
\begin{proof}
Replace $F$ and $G$ in Lemma~\ref{Fublem}, with $\abs{F(z)}$ and $\mu (z)  \abs{G(z)}$, where
\[0\le \mu (z) = \frac{\int_\gamma   \abs{\varphi}^{\nicefrac{1}{2}} \abs{F} }{\int_\gamma   \abs{\varphi}^{\nicefrac{1}{2}} \abs{G} } \le 1 \quad \mbox{for all $z\in \gamma$.} \qedhere\]
\end{proof}

Given a quadratic holomorphic differential $\varphi\, \dtext z^2$ we define two partial differential operators, called the {\it horizontal} and {\it vertical derivatives}
\begin{align*}
\partial_{_\mathsf H} = \frac{\partial}{\partial z} + \frac{\varphi}{\abs{\varphi}} \frac{\partial}{\partial \bar z} \qquad \mbox{ and } \qquad
\partial_{_\mathsf V} = \frac{\partial}{\partial z} - \frac{\varphi}{\abs{\varphi}} \frac{\partial}{\partial \bar z}.
\end{align*}
If $h$ satisfies the Hopf-Laplace equation $h_z \overline{h_{\bar z}}= \varphi$, then the horizontal and vertical trajectories of $\varphi\, \dtext z^2$ are the lines of maximal and minimal stretch for $h$. Precisely, the   following identities  hold.
\begin{align} & \abs{ \partial_{_\mathsf H} h } = \abs{h_z} + \abs{h_{\bar z}}, \qquad \abs{ \partial_{_\mathsf V} h } = \big|\abs{h_z} - \abs{h_{\bar z}}\big|  \label{hHhV} \\
& \abs{ \partial_{_\mathsf H} h } \cdot \abs{ \partial_{_\mathsf V} h }= \abs{J_h}, \qquad   \abs{ \partial_{_\mathsf H} h }^2-  \abs{ \partial_{_\mathsf V} h }^2=4 \abs{\varphi} \label{hHhVJh}
\end{align}
As a consequence
\begin{equation}\label{becareful}
\abs{ \partial_{_\mathsf V} h }^2 \le \abs{J_h} \le \abs{ \partial_{_\mathsf H} h }^2.
\end{equation}

\section{Examples}\label{examples:section}

\begin{proof}[Mappings in  Example~\ref{exnoimp}]
 Actually such  solutions  can be defined in the entire plane. First we define $h$ in the upper half plane, $ \im z > 0$, where one can settle the analytic branches of power functions and the logarithm.
   \begin{equation}\label{EX1copied}
 h(z) =
\begin{cases}
  \dfrac{z^{1-\alpha}}{1-\alpha} \; +\;\dfrac{\overline{z}^{1+\alpha}}{1+\alpha} \;, \quad &\textnormal{where} \;\;\;        \; \alpha = \dfrac{2}{p} \neq 1
  \phantom{\bigg|} \\
   \log z \;+ \dfrac{\bar{z}^2}{2} \;, \quad\quad \;\;&\textnormal{if} \;\;\; p = 2.
\end{cases}
\end{equation}
We have
\[
\begin{split}
h_z &= z^{-\frac{2}{p}}\;,\quad\quad \textnormal{which belongs to} \, \mathscr L^p_{\textnormal {weak} } (\mathbb D)\;\textnormal{but not to} \, \mathscr L^p (\mathbb D) \\
 \overline{ h_{\bar{z}}} &= {z}^{\frac{2}{p}}\;,\quad\quad \textnormal{which belongs to} \, \mathscr L^\infty (\mathbb D) \subset  \mathscr L^p (\mathbb D).
\end{split}
 \]
  Thus the Hopf-Laplace equation $h_z \overline{h_{\bar z}} \equiv 1$ holds in the upper half plane. Then we extend $h$  to the lower half of the plane by setting $h(z) = h(\overline{z})$ for  $\im z < 0$. It is a general fact, and easy to see,  that such an extension gives a Sobolev function in the entire plane. The Hopf-Laplace equation remains true in the lower half of the plane as well.

\begin{center}\begin{figure}[h]
\includegraphics[width=0.95\textwidth]{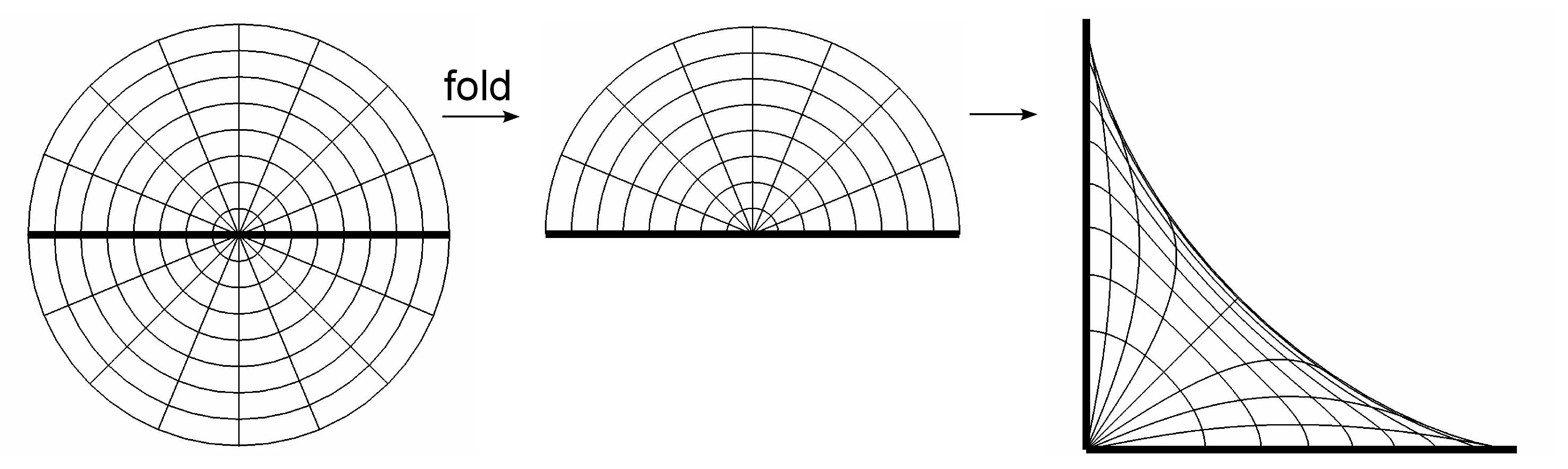} \caption{A non-Lipschitz $\mathscr W^{1,2}$-solution to the Hopf-Laplace equation}\label{fig1}
\end{figure}\end{center}

Figure~\ref{fig1} illustrates the case $p=4$. Thus $h(z)= 2z^{\nicefrac{1}{2}}+ \nicefrac{2}{3} \bar z^{\nicefrac{3}{2}}$ belongs to $\W^{1,s} (\mathbb D) \subset \mathscr W^{1,2}(\mathbb D)$ for every $2<s<4$, but is not locally Lipschitz continuous.
\end{proof}

\begin{example}\label{exsqueeze}
We use the polar coordinates for $z$ in the closed unit disk $\overline{\mathbb D}$, $z= \rho e^{i \theta}$, $0 \le \rho \le 1$ and $0 \le \theta < 2 \pi$. Define $h \colon \overline{\mathbb D} \to \C$
\[
h(\rho e^{i \theta})= 2 \rho \left[\sqrt{\rho} \sin (\nicefrac{3}{2}\, \theta) + i \sin \theta  \right]
= z- \bar z - i \left[ z^{\nicefrac{3}{2}}- \bar z^{\nicefrac{3}{2}} \right].
\]
This mapping is Lipschitz continuous, since it has bounded derivatives
\begin{equation}\label{exzbarz}
h_z = 1-\nicefrac{3}{2}\, i  \sqrt{z}, \qquad h_{\bar z} = -1+\nicefrac{3}{2}\, i  \sqrt{\bar z}.
\end{equation}
Moreover, its Hopf differential is holomorphic,
$h_z \overline{h_{\bar z}} = - \nicefrac{1}{4} \,  \left(4+  9z \right)$.
Thus $h$ solves the Hopf-Laplace equation $\frac{\partial}{\partial \bar z} \left(h_z \overline{h_{\bar z}}\right)=0$.

Formulas~\eqref{exzbarz} show that $h$  fails to be  $\mathscr C^1$-smooth in any neighborhood of  the ray $\mathbb I = \{z \colon \im z =0 \mbox{ and } 0 \le \re z \le 1\}$. Concerning topological behavior, $h$ turns out to be a harmonic diffeomorphism of $\mathbb D \setminus \mathbb I$ onto the butterfly domain $\Y \subset \C$. Figure~\ref{fig2} shows the grid of horizontal and vertical trajectories in $\X$ as well as their images in $\Y$.

\begin{center}\begin{figure}[h]
\includegraphics[width=0.8\textwidth]{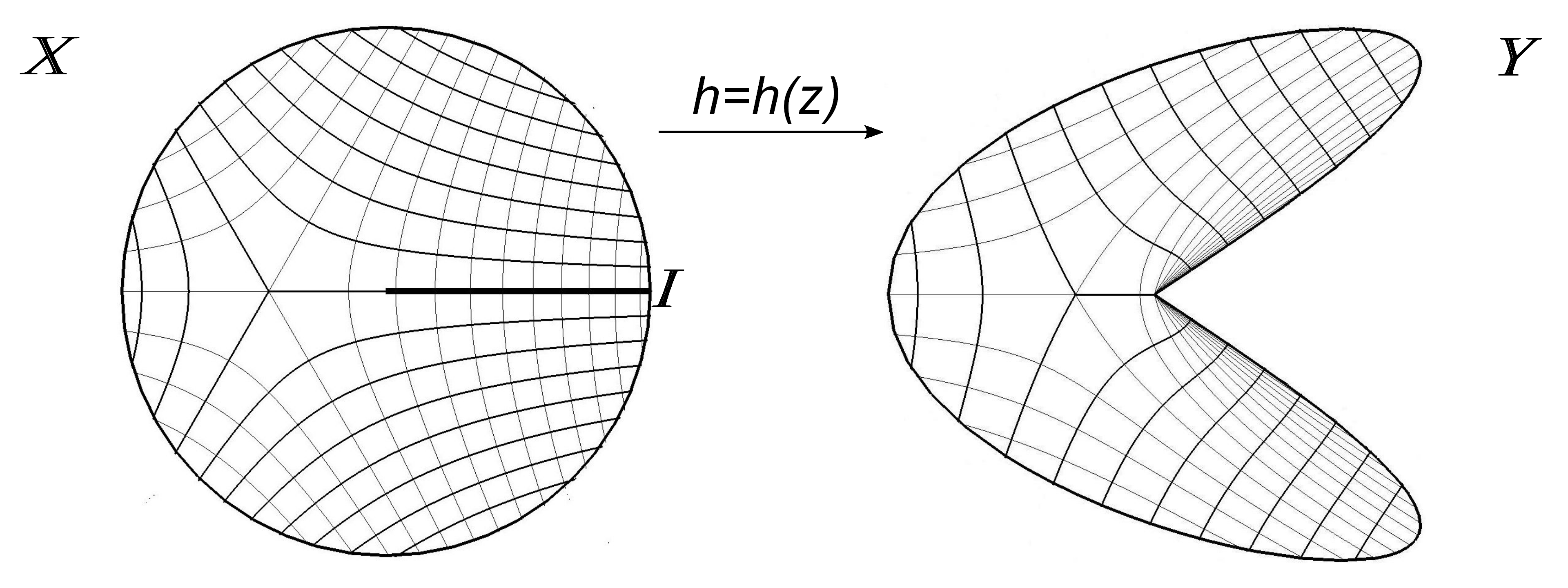} \caption{A non-$\mathscr C^1$ Hopf deformation}\label{fig2}
\end{figure}\end{center}

The radius $\mathbb I$ is squeezed  into  the origin, which is a boundary point of $\Y$. Figure~\ref{fig2} illustrates that
$\mathbb I$ is an arc of a critical vertical trajectory of the quadratic differential $\varphi\, d z^2$.
Let us notice that the functions $\abs{h_z}$ and $\abs{h_{\bar z}}$ are actually continuous. Indeed, we have
 \[\abs{h_z}^2 = 1+ \nicefrac{9}{4}\, \rho + 3 \sqrt{\rho} \sin \frac{\theta}{2}, \qquad \abs{h_{\bar z}}^2 = 1+ \nicefrac{9}{4}\, \rho - 3 \sqrt{\rho} \sin \frac{\theta}{2}.\]
In particular, the Jacobian determinant and the energy density function are also continuous
\[
\begin{split}
J_h&= \abs{h_z}^2 - \abs{h_{\bar z}}^2= 6 \sqrt{\rho} \sin \frac{\theta}{2}, \quad \mbox{ which is positive expect for $z\in \mathbb I$.}\\
\abs{Dh}^2 &= 2 \left(\abs{h_z}^2 + \abs{h_{\bar z}}^2\right)= 4 + 9 \rho.
\end{split}
\]
It is easy to see that $h$ is a $cd$-limit of homeomorphisms of $\mathbb D$ onto $\Y$. Thus $h$ is a Hopf deformation. In section~\ref{secsmo} we demonstrate through nonexplicit examples, that even minimal deformations need not be $\mathscr C^1$-smooth.
\end{example}

\section{Approximation of monotone mappings}\label{secapp}

Let us begin by recalling   the following approximation of homeomorphisms, established in~\cite{IKOhopf}.

\begin{proposition}
Let $H \colon \X \onto \Y$ be a homeomorphism of Sobolev class $\W^{1,2}_{\loc} (\X\!\shortrightarrow\!\Y)$. Then there exist diffeomorphisms $H_k \colon \X \onto \Y$, $k=1,2, \dots$, such that
\begin{itemize}
\item $H_k - H \in \mathscr A_\circ (\X)$.
\item $\norm{H_k -H}_{\mathscr A (\X)} \to 0$ as $k \to \infty$.
\end{itemize}
\end{proposition}
In view of this result we need only construct, for every $\epsilon >0$, a homeomorphism $H \colon \X \onto \Y$ such that
\begin{enumerate}[(i)]
\item $H - h \in \mathscr A_\circ (\X)$.
\item\label{ii} $\norm{H -h}_{\mathscr A (\X)} \le 6 \epsilon$.
\end{enumerate}
The construction of $H$ proceeds in three steps. We construct $\MPS(\X\!\shortrightarrow\!\Y)$ mappings $H_\circ =h$, $H_1 \in H_0 + \mathscr A_\circ (\X)$, $H_2 \in H_1 + \mathscr A_\circ (\X)$ and $H_3 \in H_2 + \mathscr A_\circ (\X)$, in which $H_3$ will turn out to be a desired homeomorphism of $\X$ onto $\Y$. In each step we make suitable  harmonic replacements to gain more points of injectivity. Moreover, estimate~\eqref{ii} will follow from:
\[\norm{H_1 -H_0}_{\mathscr A (\X)} \le 2\epsilon, \quad \norm{H_2 -H_1}_{\mathscr A (\X)} \le 2\epsilon \quad \mbox{and} \quad \norm{H_3 -H_2}_{\mathscr A (\X)} \le 2\epsilon .\]
We shall go into the construction of $H_1$ in detail in Step 1. For $H_2$ we follow the construction from Step 1, but with $H_1$ in place of $h$. In much the same way $H_3$ will be obtained as a refinement of the mapping $H_2$. Before passing to the actual construction of $H_1$ we  need some geometric considerations.
\subsection{Proof of Theorem~\ref{thm:diffeoapprox}}
An open dyadic square in $\R^2$ is the set
\[Q^m_{ij} = \{(a,b) \colon 2^m i < a < 2^m (i+1)\quad  \mbox{and} \quad   2^m j < b < 2^m (j+1)\}.\]
Hereafter, the number $2^m$ is the size of the square. Note that:
\begin{enumerate}
\item[(1)] Two different squares of the same size are disjoint.
\item[(2)] Each square of size $2^m$ is contained in exactly one square of size $2^{m+1}$, namely
\[Q^m_{ij} \subset Q^{m+1}_{\ddot \imath \,  \ddot \jmath}, \quad \mbox{ where }\;\; i-1\le 2\, {\ddot \imath} \le i \;\;\mbox{ and }\;\; j-1 \le 2 \, \ddot \jmath\le j.\]
\end{enumerate}
We refer to $Q^{m+1}_{\ddot \imath \, \ddot \jmath}$ as the dyadic square  {\it next} to $Q^m_{ij}$.
\begin{enumerate}
\item[(3)] Every two dyadic squares are either disjoint or one contains the other.
\end{enumerate}
A {\it dyadic mesh} in $\R^2$ is a family $\mathcal M$ of open dyadic squares. Let $\Y$ be a bounded domain. We will be interested only in those  dyadic squares which are compactly contained in $\Y$.  Call such a dyadic square $Q \Subset \Y$ maximal if the next dyadic square to $Q$ is not compactly contained in $\Y$. Denote by
\[\mathcal M (\Y) \quad \mbox{-the family of maximal dyadic squares in $\Y$.}\]
Clearly $\mathcal M (\Y)$ is a disjoint family and
\[\Y = \bigcup_{Q \in \mathcal M (\Y)}  \overline{Q}.\]
{\bf Claim 1.} {\it Every compact subset $\mathbb F \subset \Y$ intersects at most a finite number of closed  squares $\overline{Q}$, where $Q\in \mathcal M (\Y)$.}
\begin{proof}
For, if not, we would find an arbitrarily small square $Q\in \mathcal M (\Y)$ whose closure intersects $\mathbb F$, because one can accommodate  only a finite number of large squares in $\Y$. But then the next  dyadic square, being small enough, would be compactly contained in $\Y$. This contradicts  maximality  of $Q\in \mathcal M (\Y)$.
\end{proof}
In what follows we shall subdivide each $Q\in \mathcal M (\Y)$ into $4^n$ congruent dyadic subsquares, later referred to as fine squares. The numbers $n=n_Q$ will be chosen and fixed  according to the needs for the construction of $H_1$. In the meantime, let us reserve a notation and point out basic features of fine squares.
\[ Q_\alpha \subset Q , \qquad \alpha =1,2, \dots , 4^n, \quad n=n_Q. \]
These are disjoint open dyadic squares of  size $2^{m-n}$, where $2^m$ is the size of $Q$. For each $ Q \in \mathcal M (\Y)$ we have
\[\overline{Q}= \bigcup_{\alpha =1}^{4^n} \overline{Q}_\alpha \, , \qquad n=n_Q.\]
Once a subdivision of each  $Q\in \mathcal M (\Y)$ is made the family of fine squares will be denoted by
\[\mathcal F (\Y) = \{ Q_\alpha \colon Q \in \mathcal M (\Y), \quad Q_\alpha \subset Q, \quad \alpha =1, \dots , 4^n, \quad n=n_Q \}.\]
This is a disjoint family of  squares whose closures cover the entire domain $\Y$. As in Claim 1, we have

{\bf Claim 2.} {\it Every compact set $\mathbb F \subset \Y$ intersects at most a finite number of closed fine squares.}

We now proceed to the  construction of  the mapping $H_1$.

{\bf Step 1.} Let $Q \Subset \Y$ be a generic square in $\mathcal M (\Y)$, and $Q_\alpha \subset Q$, $\alpha =1, \dots , 4^n$, the corresponding fine  squares, with $n=n_Q$ to be determined later. Since $h\in \MPS(\X\!\shortrightarrow\!\Y)$, by Proposition~\ref{homext}, each preimage
\[\mathbb U_\alpha = h^{-1} (Q_\alpha) \Subset \X , \qquad \alpha =1,2, \dots , 4^n\]
 is a simply connected domain compactly contained in $\X$. We refer to $\mathbb U_\alpha$ as cells in $\X$. Caveat lector---the closed cell $\overline{\mathbb U}_\alpha$ can be substantially smaller than $h^{-1} (\overline{Q}_\alpha)$; it may even lie in the interior of $h^{-1} (\overline{Q}_\alpha)$.
The Youngs refinement, Proposition~\ref{homext}, tells us that $h \colon \mathbb U_\alpha \onto Q_\alpha$ admits a homeomorphism $h_{Q_\alpha} \colon  \mathbb U_\alpha \onto Q_\alpha$ with continuous extension to $\X$. The extended mapping, still denoted by  $h_{Q_\alpha} \colon  \X \onto \Y$, coincides with $h$ on $\X \setminus \U_\alpha$ and belongs to $\MPS(\X\!\shortrightarrow\!\Y)$.
However, the Youngs refinement does not guarantee that $h_{Q_\alpha}$ belongs to $\W^{1,2}_{\loc} (\X\!\shortrightarrow\!\Y)$. At this point, since  $h \in \W^{1,2}_{\loc} (\X\!\shortrightarrow\!\Y)$, the Poisson operator in Proposition~\ref{poisson} comes to the rescue. We simply replace each $h_{Q_\alpha} \colon  \mathbb U_\alpha \onto Q_\alpha$ with a harmonic diffeomorphism \[h_\alpha := {\bf P}_{_{\mathbb U_\alpha}} (h_{Q_\alpha}) \colon  \mathbb U_\alpha \onto Q_\alpha, \quad h_\alpha \in h+ \mathscr A_\circ (\mathbb U_\alpha) \]
and extend to $\X$ by setting $h_\alpha =h$ on $\X \setminus \mathbb U_\alpha$.
This yields a mapping
\[h_Q^n \colon \U \onto Q , \quad \mbox{defined by } h_Q^n = h+ \sum_{\alpha =1}^{4^n} [h_\alpha -h]_\circ\]
where $[h_\alpha -h]_\circ$ stands for the function in $\X$ that equals $h_\alpha -h$ in $\U_\alpha$ and vanishes outside $\U_\alpha$. The energy of each $h_\alpha$ does not exceed that of $h$, namely  $\E_{\U_\alpha} [h_\alpha]\le \E_{\U_\alpha} [h]$. Outside of the cells $\U_\alpha$ the mapping $h_Q^n$ coincides with $h$, thus has the same energy as $h$. Therefore,
\begin{enumerate}
\item $h_Q^n -h \in \mathscr A_\circ (\U)$.
\item $\E_{\U} [h_Q^n] \le \E_\U [h]$.
\end{enumerate}
We also have
\[\norm{h_Q^n -h}_{\mathscr C (\U)}  \le \max_{1\le \alpha \le 4^n}  \norm{h_\alpha -h}_{\mathscr C (\U)} \le \max_{1\le \alpha \le 4^n} \diam Q_\alpha = 2^{-n} \diam Q.\]
When $n$ increases to $\infty$ the mappings $h^n_Q$ converge uniformly to $h$ on $\overline{\U}$. Furthermore,  they are bounded in  $\W^{1,2} (\U)$. Thus $h^n_Q$  converge weakly to $h$ in  $\W^{1,2} (\U)$.
By weak lower semicontinuity, we have
\[\E_\U [h] \le \liminf_{n \to \infty} \E_\U [h_Q^n] \le \E_\U [h]\]
so $ \E_\U [h_Q^n]  \to  \E_\U [h] $.  We now recall the well known fact that if functions in $\mathscr L^2 (\U)$ converge weakly and their norms converge to the norm of the weak limit then such functions actually converge strongly.

 It is at this stage that we  choose and fix number $n=n_Q$, which will also depend on $\epsilon$, to be large enough to satisfy
\[
\begin{split}
\norm{h_Q^n -h}_{\mathscr C (\U)} & \le 2^{-n_Q} \diam Q \le \epsilon\\
\E_\U [h_Q^n -h] &\le \frac{\abs{Q} \, \epsilon^2 }{\abs{\Y}}.
\end{split}
\]
Finally, we conjoin all mappings $h_Q^n \colon \U \to Q$, with  $Q\in \mathcal M (\Y)$ and $n=n_Q(\epsilon)$. We obtain the desired mapping $H_1 \colon \X \onto \Y$,
\[H_1 = h+ \sum_{Q \in \mathcal M (\Y)} [h_Q^n -h]_\circ \in \mathscr A_\circ (\X) .\]
Clearly, we have
\[
\begin{split}
\norm{H_1 -h}_{\mathscr C (\X)} &=\sup_{Q\in \mathcal M (\Y)}  \norm{h_Q^n -h}_{\mathscr C (\U)} \le \epsilon \\
\E_\X [H_1-h] &= \sum_{Q\in \mathcal M (\Y)} \E_{\U} [h_Q^n -h] \le  \epsilon^2 \sum_{Q \in \mathcal M (\Y)} \frac{\abs{Q}}{\abs{\Y}} = \epsilon^2.
\end{split}
\]
Hence the estimate,
\[\norm{H_1 -h}_{\mathscr A (\X)} \le 2 \epsilon .\]
What we gained, as compared to $h$, is that the mapping $H_1 \colon \X \onto \Y$ is a harmonic diffeomorphism on every cell
\[\U_\alpha = h^{-1}_\alpha (Q_\alpha) \subset \X, \quad Q \in \mathcal M (\Y),  \quad \alpha =1,2, \dots , n_Q. \]
Thus
\begin{equation}\label{star}
H_1^{-1} (y) \quad \mbox{is a singleton if } y \in \bigcup_{Q_\alpha\in \mathcal F (\Y)} Q_\alpha.
\end{equation}
For other preimages, we have
\begin{equation}\label{sstar}
H_1^{-1} (y)=h^{-1}(y) \quad \mbox{if } y \in \bigcup_{Q_\alpha \in \mathcal F (\Y)} \partial Q_\alpha.
\end{equation}
In either case the preimage of a point in $\Y$ is connected. Thus $H_1$ is a monotone mapping. Similarly we argue that $H_1$ is a proper mapping. Indeed, let $\mathbb F$ be compact in $\Y$.  There are only finite number of closed fine  squares in  $\mathcal F(\Y)$ which  intersect $\mathbb F$. Therefore the preimage of $\mathbb F$ under $H_1$ is contained in the union of
  a finite number of closed cells in $\X$ and in  $h^{-1} (\mathbb F)$. Thus $H_1^{-1} (\mathbb F)$ stays away from $\partial \X$ and, being relatively closed in $\X$, is indeed compact.  Step~1 is completed.

 {\bf Steps 2 and 3.}
In Step 1  the construction of  $H_1 \colon \X \onto \Y$ started with a mesh $\mathcal M$ of dyadic squares in $\R^2$;  let us now  redenote this mesh as $\mathcal M_1$. Such a mesh actually depends on the choice of the orthogonal coordinates for $\R^2$. Translating the origin of the coordinate system leads to  new meshes.  These meshes would work for the constructions of $H_1$ just as $\mathcal M_1$. It would, however, lead us to different dyadic  squares in $\Y$, different family $\mathcal F (\Y)$, and different cells in $\X$. We shall take advantage of this observation by  considering  three {\it incommensurate} meshes in $\R^2$. One way to construct incommensurate meshes is by shifting the squares in  $\mathcal M$ through a vector with irrational coordinates, say ${\bf v}= (\sqrt{2}, \sqrt{2}) \in \R^2$. Specifically, let
\[
\mathcal M_1 = \mathcal M,\qquad \mathcal M_2 = \{Q+{\bf v} \colon Q \in \mathcal M  \},\qquad
\mathcal M_3 = \{Q-{\bf v} \colon Q \in \mathcal M  \}.
\]
The key observation is that no three squares from different meshes have a common boundary point; that is,
\begin{equation}\label{2star}
\partial Q^1 \cap \partial Q^2 \cap \partial Q^3 = \varnothing
\end{equation}
whenever $Q^1 \in \mathcal M_1$, $Q^2 \in \mathcal M_2$ and $Q^3 \in \mathcal M_3$.

Recall the corresponding families of open fine squares   $\mathcal F_1(\Y) \subset \mathcal M_1$, $\mathcal F_2(\Y) \subset \mathcal M_2$ and $\mathcal F_3(\Y) \subset \mathcal M_3$. In each family the  closures of fine squares  cover the entire domain $\Y$. But the essential feature of these families is that the open fine squares all together cover $\Y$, in symbols
\[
\Y= \bigcup \mathcal F_1(\Y)\,  \cup \,  \bigcup \mathcal F_2(\Y)\,  \cup \,  \bigcup \mathcal F_3(\Y).
\]
 We are now ready for the construction of $H_2$ and $H_3$. Following the construction of $H_1$ in  Step~1, but  with $H_1$ in place of $h$ and with mesh $\mathcal M_2$ in place of $\mathcal M\,,$  we obtain a mapping $H_2 \colon \X \onto \Y$ which is continuous monotone and proper. Moreover,
\[H_2-H_1 \in \mathscr A_\circ (X)  \quad \mbox{ and } \quad  \norm{H_2 -H_1}_{\mathscr A (\X)}  \le 2 \epsilon.\]
Then, in the same fashion, we refine $H_2$ by using the mesh $\mathcal M_3$.  We  arrive  at the desired  mapping $H_3 \colon \X \onto \Y$ such that $\norm{H_3-H_2}_{\mathscr A (\X)} \le 2 \epsilon$. To see that $H_3$ is a homeomorphism we need only verify injectivity. Let $y\in \Y = \bigcup_{Q\in \mathcal F_3 (\Y)} \overline{Q}$ and suppose, to the contrary, that $H_3^{-1}(y) \subset \X$ is not a singleton. This means that $y \notin \bigcup_{Q\in \mathcal F_3 (\Y)} Q$ so $y$ lies in the boundary of some square $Q^3 \in \mathcal F_3 (\Y) \subset \mathcal M_3$. Recall that for such a boundary point we have $H_3^{-1} (y) = H^{-1}_2(y)$. This in turn means that  $y \notin \bigcup_{Q\in \mathcal F_2 (\Y)} Q$, so $y$ lies in the boundary of a square $Q^2 \in  \mathcal F_2 (\Y) \subset \mathcal M_2$. For such a point we have $H_2^{-1} (y)= H_1^{-1} (y)$. As before, this means that   $y \notin \bigcup_{Q\in \mathcal F_1 (\Y)} Q$, so $y$ lies in the boundary of a square $Q^1 \in \mathcal F_1 (\Y) \subset \mathcal M_1$. In conclusion, $y$ belongs to
$\partial Q^1 \cap \partial Q^2 \cap \partial Q^3$. This contradicts~\eqref{2star}. The proof of Theorem~\ref{thm:diffeoapprox} is complete.\qed

Next we apply this theorem to harmonic $\MPS$ mappings.

\begin{proposition}\label{harmdiff}
Any harmonic mapping of class  $h \in \MPS(\X\!\shortrightarrow\!\Y)$ is a diffeomorphism.
\end{proposition}

\begin{proof}
Indeed, suppose $J_h= \abs{h_z}^2-\abs{h_{\bar z}}^2 \ge 0$, where we note that $h_z$ and $\overline{h_{\bar z}}$ are holomorphic functions. Since $h$ is surjective, $J_h \not\equiv 0$, which means that $h_z$ admits only isolated zeros. Then we obtain a meromorphic function  $\nu:=\nicefrac{\overline{h_{\bar z}}}{h_z}$  which is bounded by $1$. The zeros of $h_z$ are removable singularities. The maximum principle yields $\abs{\nu (z)} <1$ in $\X$; because it cannot be that $\abs{\nu (z)} \equiv 1$. Thus $J_h >0$ everywhere in $\X$ and, therefore, $h$ is a local diffeomorphism. The monotonicity implies that $h$ is actually injective.
\end{proof}

\subsection{Deformations of $\X$ into $\Y$}
Such deformations can now be completely characterized as  follows.

\begin{theorem}\label{thmequiv2}
Let  $f \colon \X \overset{\textnormal{\tiny{into}}}{\longrightarrow} \Y$ be a continuous mapping between $\,\ell$-connected bounded domains in the Sobolev class $\W^{1,2}(\X\!\shortrightarrow\!\Y)$ with nonnegative Jacobian. Then the following seven statements are equivalent:
\begin{itemize}
\item [\ding{172}]  $f$ is monotone proper and surjective.
\item [\ding{173}] $f$ is a uniform and strong $\,\mathscr W^{1,2}-$limit of  homeomorphisms $f_k \colon \X \onto \Y$,  in which $ f_k \in f + \mathscr W^{1,2}_\circ(\X\!\shortrightarrow\!\Y)\,$,  for all $\,k = 1,2 , ...\;$.
\item [\ding{174}] $f$ is a uniform and  strong $\,\mathscr W^{1,2}-$limit of $\,\mathscr C^\infty-$diffeomorphisms {$f_k \colon \X \onto \Y$}, in which $ f_k \in f + \mathscr W^{1,2}_\circ(\X\!\shortrightarrow\!\Y)\,$,  for  all $\,k = 1,2 , ...\;$.
\item [\ding{175}] $f$ is a uniform and  weak $\mathscr W^{1,2}-$limit of $\mathscr C^\infty-$diffeomorphisms {$f_k \colon \X \onto \Y$},  in which $ f_k \in f + \mathscr W^{1,2}_\circ(\X\!\shortrightarrow\!\Y)\,$,  for  all $\,k = 1,2 , ...\;$.
\item [\ding{176}] $f$ is a uniform and  weak $\mathscr W^{1,2}-$limit of  homeomorphisms $f_k \colon \X \onto \Y$, in which $ f_k \in f + \mathscr W^{1,2}_\circ(\X\!\shortrightarrow\!\Y)\,$,  for all $\,k = 1,2 , ...\;$.
\item [\ding{177}] $f$ is a  weak $\mathscr W^{1,2}-$limit of  homeomorphisms $f_k \colon \X \onto \Y$, in which $ f_k \in f + \mathscr W^{1,2}_\circ(\X\!\shortrightarrow\!\Y)\,$,  for all $\,k = 1,2 , ...\;$.
\item [\ding{178}]  $f$ is a deformation.
\end{itemize}
\end{theorem}

The assumption $J_f\ge 0$ a.e. does not impose an essential restriction on the mapping $f$, by virtue of Corollary~\ref{jacdncs}.

\begin{proof}
The implications \ding{172}~$\Rightarrow$~\ding{173}~$\Rightarrow$~\ding{174} are just a restatement of Theorem~\ref{thm:diffeoapprox}. The implications  \ding{174}~$\Rightarrow$~\ding{175}~$\Rightarrow$~\ding{176}~$\Rightarrow$~\ding{177} are obvious.

For the proof of \ding{177}~$\Rightarrow$~\ding{178} we argue as follows.  Since $f_k \colon \X \onto \Y$ are orientation preserving homeomorphisms of Sobolev class $f + \W^{1,2}_\circ (\X\!\shortrightarrow\!\Y)$ converging weakly in $\W^{1,2}(\X\!\shortrightarrow\!\Y)$ to $f$,
it follows that
\[\iint_\X J_f \, dx =  \iint_\X J_{f_k} \, dx= \abs{\Y}. \]
Hence there is a compact subset $\mathbb F \subset \X$ such that
\[\iint_{\mathbb F} J_f \, dx > \frac{1}{2} \abs{\Y}.\]
By weak $\mathscr L^1$-convergence of nonnegative Jacobians~\cite[Theorem 8.4.2]{IMb} we have $\iint_{\mathbb F} J_f \, dx = \lim \iint_{\mathbb F} J_{f_k}\, dx$. Therefore,
\[\abs{f_k (\mathbb F)} = \iint_{\mathbb F} J_{f_k}\, dx > \frac{1}{2} \abs{\Y}, \quad \mbox{for sufficiently large $k$.}\]
It then follows that there is $\epsilon >0$, such that
\[\sup_{x\in \mathbb F} \left(f_k(x), \partial \Y\right) \ge \epsilon, \qquad k=1,2, \dots .\]
Let $x_k \in \mathbb F$ be a point for which $y_k=f_k (x_k)$ has distance at least $\epsilon$ from $\partial \Y$. We may assume, by passing to a subsequence if necessary, that
\[x_k \to x_\circ \in \mathbb F \qquad \mbox{and} \qquad y_k \to y_\circ \in \Y.\]
Let $\Phi_k \colon \X \onto \X$ and $\Psi_k \colon \Y \onto \Y$ be local perturbations (arbitrarily small) of the identity mapping near $x_\circ$ and $y_\circ$, respectively, to satisfy,
\[\Phi_k (x_\circ)=x_k  \qquad \mbox{and} \qquad  \Psi_k (y_k)=y_\circ. \]
Now the homeomorphisms $F_k = \Psi_k \circ f_k \circ \Phi_k \colon \X \onto \Y$ take $x_\circ$ into $y_\circ$. We consider $F_k$ as a mapping of the punctured domain $\X_\circ = \X \setminus \{x_\circ\}$ onto $\Y_\circ = \Y \setminus \{y_\circ\}$, each of which has $\ell +1 \ge 2$ boundary components. These mappings coincide with $f_k$ outside a compact subset of  $\X$. At this point we appeal to the following uniform estimate of the distance to $\partial \Y_\circ$~\cite[Theorem 1.1]{IOt}.
\begin{equation}\label{distest}
\dist \left(F_k (x), \partial \Y_\circ \right) \le \eta (x)\, \norm{DF_k}_{\mathscr L^2 (\X_\circ)}, \qquad k=1,2, \cdots
\end{equation}
where $\eta (x)= \eta_{_{\X_\circ \, \Y_\circ}} (x)$ is a continuous function in $\overline{\X_\circ}$ vanishing on $\partial \X_\circ$. We emphasize that this function depends only on the domains $\X_\circ$ and $\Y_\circ$. Since $F_k = f_k$ near $\partial \X$ the estimate~\eqref{distest} yields, for each $f_k$,
\[\dist \left(f_k (x), \partial \Y\right) \le \eta (x)\, M\]
where $M$ is controlled from above by the energy of $f_k$, so is independent of $k$. Finally, since $f_k \to f$ $c$-uniformly we conclude that $\dist (f_k (x), \partial \Y) \to \dist (f (x), \partial \Y) $ uniformly in $\X$. This shows that $f$ is a deformation.

The implication \ding{178}~$\Rightarrow$~\ding{172} is a part of Proposition~\ref{uniappthm}.
\end{proof}

\begin{remark}
The observant reader may notice that the conditions \ding{173}--\ding{177} tell us something about the boundary behavior of a deformation with range $\Y$. In a way every deformation $f \colon \X \onto \Y$ must agree on $\partial \X$ with a homeomorphism $f_k \colon \X \onto \Y$. This is understood in the sense of Sobolev boundary data  $f\in f_k + \W^{1,2}_\circ (\X\!\shortrightarrow\!\Y)$.
\end{remark}

\section{Preimage of a point under generalized solutions}\label{secpreimages}

Let $h \colon \X \to \R^2$ be a continuous mapping. The multiplicity function of $h$ defined by $\mathcal N_h(y)= \# \{h^{-1} (y)\}$, $y\in \R^2$, is measurable, so one can speak of the essential supremum of $\mathcal N_h (y)$. We are concerned with mappings such that
\begin{equation}\label{muless} \underset{y\in \R^2}{\esssup} \, \mathcal N_h (y)< \infty .\end{equation}
Note that Hopf deformations enjoy the property~\cite[Lemma 3.8]{IKKO}
\begin{equation}\label{enjoyable}
\underset{y\in \R^2}{\esssup} \, \mathcal N_h (y)= 1.
\end{equation}
The following proposition deals with more general solutions to the Hopf equation.

\begin{proposition}\label{lily}
Let $h \colon \X \to \C$ be a continuous $\W^{1,1}_{\loc} (\X)$-solution to the Hopf-Laplace equation
\[h_z \overline{h_{\bar z}}= \varphi \not\equiv 0 \qquad \mbox{almost everywhere in } \X,\]
where $\varphi$ is a holomorphic function in a domain $\X \subset \C$. Assume that the multiplicity function $\mathcal N_h(y)= \# \{h^{-1} (y)\}$ is essentially bounded~\eqref{muless}.
Then for each $y_\circ \in \R^2$ the union of all vertical trajectories of the quadratic differential $\varphi (z)\, \dtext z^2$ in $\X$ that intersect $h^{-1} (y_\circ)$ has zero measure.
\end{proposition}

\begin{proof}
To simplify writing we assume that $y_\circ = 0 \in \R^2$.
Let $\mathcal V$ denote the family of all vertical trajectories of $\varphi (z)\, \dtext z^2$ in $\X$. These are disjoint open $\mathscr C^\infty$-smooth curves without self-intersections whose union covers $\X_\circ = \X \setminus \{\mbox{zeros of } \varphi\}$. Every point in $\X_\circ$ has a neighborhood in which a single valued branch of the analytic function $\Phi (z)= \int \sqrt{\varphi (z)}\, \dtext z$ can be chosen. This is a local conformal mapping which takes the arcs of vertical trajectories into open vertical intervals in the $w$-plane, $w=\Phi (z)$. In general it may not be possible to perform analytic continuation of $\Phi$ along the entire trajectory; the local branches of $\int \sqrt{\varphi}$  may not coincide if  their domains of definition are overlapping. This difficulty is usually overcome by performing analytic continuation of the inverse $\Phi^{-1}$ along the straight vertical lines in the $w$-plane, see~\cite[\S1.3.2]{Stb} for a thorough discussion. Such a procedure leads to the concept of a vertical strip. A vertical strip in the $w$-plane  associated with $\varphi \, \dtext z^2$ is a simply connected domain  of the form
\[\mathcal S= \{w=t+i \tau \colon 0<t<T, \quad \alpha (t)< \tau < \beta (t)\}  \]
where $-\infty \le \alpha (t)< \beta (t) \le \infty$ are measurable functions in $t\in (0,T)$. Moreover, there is a single valued analytic function $\Psi \colon \mathcal S \into \X$ which takes every vertical interval $\gamma_t= \{t+i\tau \colon \alpha (t)< \tau < \beta (t)\}$ onto a complete vertical trajectory in $\X$. This mapping $\Psi$ is locally conformal and its  inverse, locally defined, is a branch of $\Phi = \int \sqrt{\varphi (z)}\, \dtext z$. Thus the image $\Psi (\mathcal S) \subset \X$ is an open subset of $\X$.  Each vertical trajectory in $\X$ either lies entirely in $\Psi (\mathcal S)$ or otherwise is disjoint from $\Psi (\mathcal S)$. The point is that the whole domain $\X_\circ$
 can be covered by a countable number of domains such as $\Psi (\mathcal S)$.

Denote $\mathcal V_\circ \subset \mathcal V$ the family of vertical trajectories in $\X$ which intersect the set $h^{-1}(0)$ and assume, to derive a contradiction, that the union $\bigcup \mathcal V_\circ$ has positive measure. We shall confine ourselves to one particular subdomain $\Psi (\mathcal S) \subset \X$ and trajectories selected from $\mathcal V_\circ$ that lie in $\Psi (\mathcal S)$. With a suitable choice of $\Psi (\mathcal S)$ we ensure that the union of the selected trajectories still has positive measure. Rather than discuss this subdomain, let us assume that $\X= \Psi (\mathcal S)$. Further simplification comes by considering the mapping $f= h \circ \Psi \colon \mathcal S \to \C$. This simplifies not only the domain of definition but also the Hopf-Laplace equation translates into the somewhat easier form
\begin{equation}\label{HL345} f_w \overline{f_{\bar w}} \equiv 1, \qquad \mbox{for all } w= t+i \tau \in \mathcal S.\end{equation}
Let $\Gamma= \{\gamma_t\}_{0<t<T}$ denote the family of all vertical intervals in $\mathcal S$; these are vertical trajectories of $f\, \dtext w^2$,
\[\gamma_t = \{t+i\tau \colon \alpha (t)< \tau < \beta (t)\}, \qquad 0<t<T.\]
Among them there are intervals that pass through the set $f^{-1}(0)$ which we  designate by
\[\Gamma_\circ =\{\gamma_t \in \Gamma \colon 0 \in f(\gamma_t)\}.\]
In this way we are reduced, equivalently, to showing that the union $\bigcup \Gamma_\circ$ has positive measure. That this is indeed an equivalent problem  follows from the observation that $\Psi$, being a local diffeomorphism, takes a null family of vertical intervals in $\mathcal S$ into a null family of vertical trajectories in $\Psi (\mathcal S)$. Furthermore, the strip $\mathcal S$, possible infinite, can be exhausted with an increasing sequence of vertical strips compactly contained in $\mathcal S$,
\[\mathcal S_1 \Subset \mathcal S_2 \Subset \cdots \mathcal S_n \Subset \dots \Subset \mathcal S = \bigcup_{n=1}^\infty \mathcal S_n.\]
Let the family $\Gamma_\circ^n$ consist of those vertical intervals in $\mathcal S_n$ which pass through the set $f^{-1}(0)$. Clearly, we have $\bigcup_{n=1}^\infty \left(\bigcup \Gamma_\circ^n\right)= \bigcup \Gamma_\circ$; the latter is a subset of $\mathcal S$ with positive measure. Thus for some large $n$  we still  have
\[\left| \bigcup \Gamma_\circ^n \right|>0.\]
Therefore, we may and do assume, instead of introducing new notation that $\mathcal S$ is bounded and $\Psi \colon \mathcal S \onto \Psi (\mathcal S)$ extends as a local conformal mapping to a neighborhood of $\overline{\mathcal S}$. In particular, $\abs{\mathcal S}< \infty$ and the multiplicity function of $\Psi \colon \mathcal S \to \Psi (\mathcal S)$ is also bounded.

\begin{lemma}
For almost every $t\in (0, T)$ such that $\gamma_t \in \Gamma_\circ$, we have
\[\diam f(\gamma_t)>0.\]
\end{lemma}
\begin{proof}
Let $C \subset (0,T)$ denote the set of parameters $t$ such that $\diam f(\gamma_t)=0$ and $\gamma_t \in \Gamma_\circ$. This means that $f$ is a constant mapping on each interval $\gamma_t$, for $t\in C$. Since $0\in f(\gamma_t)$ we conclude that $f\equiv 0$ on $\bigcup_{t\in C} \gamma_t$. On the other hand, in view of the Hopf-Laplace equation~\eqref{HL345}, $f$ cannot vanish on a set of positive measure, so $\left|\bigcup_{t\in C} \gamma_t\right|=0$. Hence $C$ has zero linear measure, as claimed.
\end{proof}
We now choose and denote by $E \subset (0,T)$ a set of positive linear measure such that
\begin{equation}\label{diam67}
\diam f(\gamma_t) \ge 2\rho , \qquad \mbox{for all } t \in E
\end{equation}
where $\rho$ is a sufficiently small positive number. Consider a sequence of concentric annuli centered at $0$,
\[\mathbb A_m = \{y \colon 2^{-m} \rho \le \abs{y} \le 2^{1-m} \rho  \}, \quad m=1,2, \dots\]
 It follows from~\eqref{diam67} that $f(\gamma_t)$, with  $t\in E$, is a connected set which joins $0$ with a point outside the outer boundary  of $\mathbb A_m$. Elementary geometric arguments give an  estimate of $1$-dimensional Hausdorff measure of the set  $\mathbb A_m \cap f(\gamma_t) \subset \C$, namely
\[\mathcal H^1 \left(f(\gamma_t) \cap \mathbb A_m  \right)\ge 2^{-m} \rho.\]
For almost every $t\in E$, the function $f$ is absolutely continuous on $\gamma_t$,  because $f\in \W^{1,1}(\mathcal S)$. Consider a subset $\mathbb K = \gamma_t \cap f^{-1} (\mathbb A_m)$ of the interval $\gamma_t  \subset \Gamma_\circ$. We have
\[\int_{\mathbb K} \left| \frac{\partial f}{\partial \tau}  \right|   \ge \mathcal H^1\big( f({\mathbb K})\big) \ge \mathcal H^1 \left( f(\gamma_t) \cap \mathbb A_m  \right) \ge 2^{-m} \rho \]
where we used the inclusion  $f(\mathbb K) \supset f(\gamma_t) \cap \mathbb A_m$.
Integrating with respect to $t\in E$, by Fubini's theorem, we obtain
\[\iint_{H^{-1}(\mathbb A_m)}   \left| \frac{\partial f}{\partial \tau}  \right| \ge \int_{\mathbb E} \left( \int_{\gamma_t \cap f^{-1}(\mathbb A_m)} \left| \frac{\partial f}{\partial \tau}  \right|  \right) \ge 2^{-m} \rho \, \abs{E}.\]
Next we apply H\"older's inequality
\begin{equation}\label{sub65}
4^{-m} \rho^2 \abs{E}^2 \le \abs{f^{-1} (\mathbb A_m)} \iint_{H^{-1} (\mathbb A_m)}  \left| \frac{\partial f}{\partial \tau}  \right|^2.
\end{equation}
It is at this point that we shall appeal  to the Hopf-Laplace equation~\eqref{HL345} and  formula~\eqref{becareful}, which gives us a pointwise inequality in terms of the Jacobian determinant of $f$,
\[ \left| \frac{\partial f}{\partial \tau}  \right|^2 \le \left|J_f  \right| \qquad \mbox{a.e. in }\mathcal S . \]

Now recall that the multiplicity function of $h$ is essentially bounded and $\Psi \colon \mathcal S \to \Psi (\mathcal S)$ has finite multiplicity. Therefore, the  function $\mathcal N_f (y)= \# \{w \in \mathcal S \colon f(w)=y\}$ is essentially bounded as well, say $\mathcal N_f (y) \le \mathcal N$ for almost every $y\in \R^2$.
We have
\[\iint_{f^{-1}(\mathbb A_m)}  \left| \frac{\partial f}{\partial \tau}  \right|^2 \le  \iint_{f^{-1}(\mathbb A_m)}  \abs{J_f} \le  \mathcal N\,  \abs {\mathbb A_m} = 3 \pi \rho^2 4^{-m } \mathcal N \]
where the second inequality follows from~\cite[Theorem 6.3.2]{IMb}.
Substituting into~\eqref{sub65} yields
$\abs{E}^2 \le 3 \pi \mathcal N\,  \abs{f^{-1}(\mathbb A_m)}$.
Finally we add these inequalities for $m=1,2, \dots , \ell$ to obtain
\[\ell \,  \abs{E}^2 \le 3 \pi \mathcal N\,  \Abs{f^{-1} \Big(\bigcup_{m=1}^\ell \mathbb A_m \Big)  } \le 3 \pi \mathcal N\,  \abs{\mathcal S}\]
where $\ell$ can be any positive number we wish. Thus $\abs{E}=0$, completing the proof of Proposition~\ref{lily}.
\end{proof}

\begin{corollary}\label{c1lily}
Under the assumptions of Proposition~\ref{lily}, suppose $w\in \Y$ and $h^{-1} (w)$ does not contain any critical points of $\varphi(z)\,\dtext z^2$.  Then $h^{-1} (w)$ is a closed vertical arc.
\end{corollary}
\begin{proof}
By Lemma~3.7 in~\cite{IKKO} the set $h^{-1} (w)\subset \X$ is connected and compact. Then the union of vertical trajectories which intersect $h^{-1}(w)$ is connected and, by Proposition~\ref{lily}, has zero measure. This is possible only when $h^{-1}(w)$ is contained in exactly one vertical trajectory.
\end{proof}

\section{Partial harmonicity, proof of Theorem~\ref{thharm}}

The outline of the proof is as follows. We may assume that the Hopf differential of $h$ does not vanish identically, for otherwise $h$ is holomorphic in $\X$.
Using the notation of Theorem~\ref{thharm} let $\mathbb D$ be an open disk compactly contained in $\Y$
and $\U=h^{-1} (\mathbb D)$. We will prove  that $h$ is harmonic in $\U$ by showing that the energy of $h$ does not exceed the energy of $H$,
\begin{equation}\label{goal}
\E_\U [h] \le \E_{\U} [H]
\end{equation}
where $H$ is the Poisson refinement of $h$ in $\U$.
Indeed,~\eqref{goal} shows that $h=H$ in $\U$ and therefore $h$ is harmonic diffeomorphism of $\U$ onto $\mathbb D$. Since $h$ is also monotone,
it is a global diffeomorphism. Thus~\eqref{goal} is all we need to prove Theorem~\ref{thharm}.

The first step toward proving~\eqref{goal} is the following  computation.

\begin{lemma}\label{energylemma}
Let $\U$ and $\mathbb D$ be bounded simply connected domains. Suppose that $h \colon \U \onto \mathbb D$, of Sobolev class $\mathscr W^{1,2} (\U, \mathbb D)$, is monotone and proper and has a continuous extension to $\overline{\U}$. Furthermore, let $H \colon \U \onto \mathbb D$ be a $\mathscr C^\infty$-diffeomorphism of Sobolev class $\mathscr W^{1,2} (\U)$ which extends continuously to $\overline{\U}$ with $H(z)=h(z)$ for $z\in \partial \U$.
Then  for  $\chi = H^{-1} \circ h \colon \U \onto \mathbb D$ and  $\varphi (z)= h_z \overline{h_{\bar z}}$ we have
\begin{equation}\label{compine}
\begin{split}
\E_\U [H]- \E_\U [h] &\ge \frac{4}{\norm{\varphi}_{\mathscr L^1(\U)}}\,  \left[\iint_\U \Abs{\chi_z- \frac{\varphi}{\abs{\varphi}} \chi_{\bar z}} \sqrt{\abs{\varphi(z)}}\sqrt{\abs{\varphi \big(\chi (z) \big)}} \, dz   \right]^2
\\
& - 4 \iint_{\U} \abs{\varphi}.
\end{split}
\end{equation}
Here we assume that $\varphi$ is continuous, $\varphi\not\equiv 0$, and the term $\frac{\varphi}{ \abs{\varphi}}$ is understood as equal to zero
whenever $\varphi$ vanishes.
\end{lemma}
\begin{proof}
First assume, in addition to the above hypotheses that $h \colon \U \onto \mathbb D$ is a diffeomorphism. 

The chain rule can be applied to the composition $H=h\circ \chi^{-1}\colon \U \onto \U$
\begin{equation*}
\begin{split}
\frac{\partial H(w)}{\partial w} &= h_z(z)\frac{\partial \chi^{-1}}{\partial w}+h_{\bar z}(z)  \frac{\overline{\partial \chi^{-1}} }{\partial \overline{w}} \\
\frac{\partial H(w)}{\partial \bar w} &= h_z(z)\frac{\partial \chi^{-1}}{\partial \bar w}+h_{\bar z}(z)\overline{\frac{\partial \chi^{-1}}{\partial w}}
\end{split}
\end{equation*}
where $w=\chi(z)$. The partial derivatives of $\chi^{-1} \colon \U \to \U$ at $w$ can be expressed  in terms $\chi_z$ and $\chi_{\bar z}$ at $z= \chi^{-1} (w)$ by the rules
\begin{equation*}
\frac{\partial \chi^{-1}}{\partial w} = \frac{\overline{ \chi_z (z)}}{J(z,\chi)} \quad \mbox{ and } \quad
\frac{\partial \chi^{-1}}{\partial \bar w} = - \frac{\chi_{\bar z}(z)}{J(z,\chi)}
\end{equation*}
where the Jacobian determinant $J(z,\chi)$ is strictly positive. This yields
\begin{equation*}
\frac{\partial H}{\partial w}  =\frac{h_z\overline{\chi_z}-h_{\bar z} \overline{\chi_{\bar z}} }{J(z,\chi)} \quad \mbox{ and } \quad \frac{\partial H}{\partial \bar w}  =\frac{h_{\bar z}{\chi_z}-h_{ z}{\chi_{\bar z}}}{J(z,\chi)}.
\end{equation*}
Let $\U' \Subset \U$ be a compactly contained subdomain of $\U$. We compute the energy of $H$ over the set $\chi (\U')$ by substitution $w= \chi (z)$,
\begin{equation*}
\begin{split}
\E_\U [H] \ge\mathcal E_{\chi(\U')}[H]&=2\iint_{\chi(\U')} \left(\abs{H_w}^2+\abs{H_{\bar w}}^2\right)\,\dtext w \\
&=2\iint_{\U'}\frac{\abs{h_z \overline{\chi_z} - h_{\bar z} \overline{\chi_{\bar z}}}^2+
\abs{h_{\bar z} \chi_z - h_{z} \chi_{\bar z}}^2}{\abs{\chi_z}^2-\abs{\chi_{\bar z}}^2}\,\dtext z.
\end{split}
\end{equation*}
On the other hand, the energy of $h$ over the set $\U'$ is
\[
\mathcal E_{\U'}[h]=2\iint_{\U'}\left(\abs{h_z}^2+\abs{h_{\bar z}}^2\right)\,\dtext z.
\]
Subtract these two integral expressions to obtain
\begin{equation}\label{cchain}
\begin{split}
\mathcal E_{\U}[H]-\mathcal E_{\U'}[h] &\ge  4\iint_{\U'}  \frac{\left(\abs{h_z}^2+\abs{h_{\bar z}}^2\right)\cdot \abs{\chi_{\bar z}}^2
-2\re \left[h_z\overline{h_{\bar z}} \overline{\chi_z}\chi_{\bar z}\right]}{\abs{\chi_z}^2-\abs{\chi_{\bar z}}^2}\,\dtext z \\
& \ge 4 \iint_{\U'} \frac{ 2\abs{h_zh_{\bar z}} \cdot \abs{\chi_{\bar z}}^2
-2\re \left[h_z\overline{h_{\bar z}} \overline{\chi_z}\chi_{\bar z}\right]}{\abs{\chi_z}^2-\abs{\chi_{\bar z}}^2}\,\dtext z \\
& = 4 \iint_{\U'} \left[\frac{\abs{\chi_z-\sigma(z) \chi_{\bar z}}^2}{\abs{\chi_z}^2-\abs{\chi_{\bar z}}^2} -1 \right]
\, \abs{h_zh_{\bar z}} \,\dtext z
\end{split}
\end{equation}
where
\[
\sigma = \sigma(z) = \begin{cases}
{h_z\overline{h_{\bar z}}}{\, \abs{h_z  \overline{h_{\bar z}}}^{-1}} \qquad &\text{if } h_z\overline{h_{\bar z}}\ne 0 \\
0 & \text{otherwise.}
 \end{cases}
\]
Using H\"older's inequality we continue the above chain of estimates as follows
\begin{equation}\label{cchain2}
\ge \; 4\frac{\left[\iint_{\U'}\left| \chi_z-\sigma\chi_{\bar z}\right| \,\sqrt{\abs{h_zh_{\bar z}}} \sqrt{\abs{\psi \big(\chi (z) \big)}}   \,\dtext z\right]^2}{\iint_{\U'} J(z,\chi) \abs{\psi \big( \chi(z) \big)}  \,\dtext z}
-4\iint_{\U'} \abs{h_zh_{\bar z}}.
\end{equation}
where $\psi \colon \U \to \C$ can be any continuous function, provided $\psi \not\equiv 0$ on $\chi (\U')$.

The denominator in~\eqref{cchain2} is uniformly bounded from above
\[\iint_{\U'} J(z,\chi) \abs{\psi \big( \chi(z) \big)}  \,\dtext z = \iint_{\chi (\U')} \abs{\psi} \le \iint_\U \abs{\psi}.\]
Hence
\begin{equation*}
\mathcal E_{\U}[H]-\mathcal E_{\U'}[h]\ge
4\, \frac{\left[\iint_{\U'}\left| \chi_z-\sigma\chi_{\bar z}\right| \,\sqrt{\abs{\varphi (z) }} \sqrt{\abs{\psi\big(\chi (z) \big)}}   \,\dtext z\right]^2} {\iint_{\U} \abs{\psi (z)}  \,\dtext z}
 -4\iint_{\U'}\abs{\varphi}.
\end{equation*}
 This inequality can now be generalized by an approximation argument. By Theorem~\ref{thm:diffeoapprox}, we have a sequence of diffeomorphisms
 $h^j \colon \U \onto \mathbb D$, converging to $h$ uniformly and strongly in $\mathscr W^{1,2} (\U, \mathbb D)$. Moreover, each
 $h^j \in h+ \mathscr A_\circ (\U, \mathbb D)$, so $h^j$ extends continuously to $\overline{\U}$ with $h^j(z)=h(z)$ on $\partial \U$. We may and do assume,
 by passing to a subsequence if necessary, that $h^j_z$ and $h^j_{\bar z}$ converge almost everywhere to $h_z$ and $h_{\bar z}$, respectively.
Let $\varphi^j=h_z^j \overline{h_{\bar z}^j}$.
Since the sequence $\chi^j= H^{-1} \circ h^j \colon \U \onto \U$ of self-diffeomorphisms of $\U$ is converging to $\chi$ uniformly and strongly
in $\mathscr W^{1,2}$ on subdomains $\U' \Subset \U$, it follows that
\[
\mathcal E_{\U}[H]-\mathcal E_{\U'}[h^j]\ge
4\frac{\left[\iint_{\U'}\left| \chi^j_z-\sigma^j \chi^j_{\bar z}\right| \,\sqrt{\abs{\varphi^j (z) }} \sqrt{\abs{\psi\big(\chi (z) \big)}}   \,\dtext z\right]^2} {\iint_{\U} \abs{\psi (z)}  \,\dtext z} -4\iint_{\U'}\abs{\varphi^j}.
\]
Passing to the limit as  $j \to \infty$ yields
\begin{equation}\label{eq2222}
\begin{split}
\mathcal E_{\U}[H]-\mathcal E_{\U'}[h]&\ge
4\frac{\left[\iint_{\U'}\left| \chi_z-\sigma\chi_{\bar z}\right| \,\sqrt{\abs{\varphi (z) }} \sqrt{\abs{\psi\big(\chi (z) \big)}}   \,\dtext z\right]^2} {\iint_{\U} \abs{\psi (z)}  \,\dtext z}\\ & -4\iint_{\U'}\abs{\varphi}.
\end{split}
\end{equation}
To see this we simply note that
\[ \left| \chi_z^j-\sigma^j\chi_{\bar z}^j\right| \,\sqrt{\abs{h_z^jh_{\bar z}^j}}\; \to \;  \left| \chi_z-\sigma\chi_{\bar z}\right| \,\sqrt{\abs{\varphi (z)}} \quad \mbox{ in } \mathscr L^1 (\U') \]
while $\sqrt{\abs{\psi \big(\chi^j (z)  \big)}} \to \sqrt{\abs{\psi \big(\chi (z) \big)}}$  everywhere.

Here we recall that $\varphi$ is assumed to be continuous, so we can take $\psi= \varphi$ in~\eqref{eq2222}. Finally, since $\U'$ was an arbitrary compact subset of $\U$, we conclude from~\eqref{eq2222} with the desired estimate,  completing the proof of Lemma~\ref{energylemma}.
\end{proof}

For the proof of~\eqref{goal} it remains to show that the right hand side of~\eqref{compine} is nonnegative. This requires a careful analysis of the boundary behavior of $\chi$, as
this mapping is not necessarily continuous up to $\partial\mathbb U$. We need a definition and a lemma.

\begin{definition}\label{needit}
Let $\X\subset \C$ be a domain and $\mathbb U$ be a simply connected domain compactly contained in $\X$.
Let $\varphi \colon \X \to \C$ be a holomorphic function such that $\varphi\not\equiv 0$. We say that a mapping $\chi \colon \mathbb U \to \mathbb U$ is \emph{compatible with $\varphi$} if the following holds for any
vertical arc $\gamma$ of $\varphi\,\dtext z^2$ that intersects $\mathbb U$ and has endpoints in $\X\setminus \mathbb U$.
Let $\gamma_\circ$ be a maximal subarc of $\gamma$ contained in $\mathbb U$, and denote its endpoints by $a$ and $b$. The connected components of $\gamma\setminus\gamma_\circ$ are
naturally denoted as $\gamma_a$ and $\gamma_b$. The condition we impose on $\chi$ is
\[
\chi\{a\}\subset \gamma_a\quad \text{ and } \quad \chi\{b\}\subset \gamma_b
\]
where $\chi\{a\}$ and $\chi\{b\}$ are cluster sets~\cite{CL}.
\end{definition}

\begin{lemma}\label{compatible}
The mapping $\chi$ in Lemma~\ref{energylemma} is compatible with $\varphi$ provided that $h^{-1}(\partial \mathbb D)$ contains no zeros of $\varphi$.
\end{lemma}

\begin{proof} We use the notation of Definition~\ref{needit}. According to Corollary~\ref{c1lily}, the sets $h^{-1}(h(a))$ and $h^{-1}(h(b))$ are vertical arcs,
hence subarcs of $\gamma$. The continuity of $h$ and $H$ implies that $\chi\{a\}$ is a subset of $h^{-1}(h(a))$, and similarly for $\chi\{b\}$. The claim follows.
\end{proof}

The restriction concerning zeros of $\varphi$ in Lemma~\ref{compatible} is easily fulfilled by choosing a generic radius for the disk $\mathbb D$.
We are finally ready to handle the expression in~\eqref{compine}, thus completing the proof of Theorem~\ref{thharm}. The following result
is related to the Reich-Strebel inequality~\cite{RS}, see also~\cite{MM} for a recent extension. However, in our Lemma~\ref{reich-strebel} the assumptions
on $\chi$ are different from those in~\cite{RS,MM}.

\begin{lemma}\label{reich-strebel}
Let $\X\subset \C$ be a domain and $\mathbb U$  a simply connected domain compactly contained in $\X$.
Let $\varphi \colon \X \to \C$ be a holomorphic function such that $\iint_{\X} \abs{\varphi} < \infty$.
Suppose that $\chi \in \mathscr W^{1,2}_{\loc} (\U, \U)$  is continuous  proper and compatible with $\varphi$. Then
\begin{equation}\label{eq67}
\iint_{\U} \Abs{\chi_z- \frac{\varphi}{\abs{\varphi}} \chi_{\bar z}  } \abs{\varphi}^{\nicefrac{1}{2}} \abs{\varphi \circ \chi}^{\nicefrac{1}{2}} \ge \iint_{\U} \abs{\varphi}.
\end{equation}
\end{lemma}

\begin{proof}
For almost every vertical noncritical trajectory $\gamma$ the mapping $\chi$ is locally absolutely continuous on $\gamma$. Let $\gamma_\circ$ be a maximal subarc of $\gamma$ in $\mathbb U$. Denote by $a$ and $b$ the endpoints of $\gamma_\circ$.  By the compatibility condition the curve $\beta= \chi \circ \gamma_\circ$ connects two different components of $\gamma \setminus \gamma_\circ$.  By Lemma~\ref{lengthineq} we have
\begin{equation}
\int_{\gamma_\circ} \abs{\varphi}^{\nicefrac{1}{2}} \le \int_{\beta } \abs{\varphi}^{\nicefrac{1}{2}} = \int_{\gamma_\circ} \Abs{\chi_z- \frac{\varphi}{\abs{\varphi}} \chi_{\bar z}  }  \abs{\varphi \circ \chi}^{\nicefrac{1}{2}}
\end{equation}
because $ \Abs{\chi_z- \frac{\varphi}{\abs{\varphi}} \chi_{\bar z}  } $ is the magnitude of directional derivative of $\chi$ along $\gamma$. In view of Corollary~\ref{fubcor} inequality~\eqref{eq67} follows.
\end{proof}

Combining~\eqref{eq67} with~\eqref{compine}
yields~\eqref{goal}, completing the proof of Theorem~\ref{thharm}.

\section{Lipschitz continuity, proof of Theorem~\ref{MainTh}}\label{sectionlip}

In this section  $\X$ and $\Y$ are bounded $\ell$-connected domains. Suppose $h\in \mathfrak D(\X, \Y)$ is a Hopf deformation, that is,
\begin{equation}\label{hopf1}
\varphi:=h_z \overline{h_{\bar z}}
\end{equation}
is a holomorphic function in $\X$. We shall actually prove the following explicit bound
\begin{equation}\label{addstar}
\abs{Dh(a)} \le 72\,\frac{\norm{Dh}_{\mathscr L^2 (\X)}}{\dist (a, \partial \X)} \quad \mbox{for almost every $a\in \X$.}
\end{equation}
Let $\norm{\varphi}=\iint_\X \abs{\varphi}$ and note the
pointwise inequality
\begin{equation}\label{forphi}
\abs{\varphi(z)}\le \frac{\norm{\varphi}}{\pi \dist^2(z,\partial \X)}
\end{equation}
which is a consequence of the subharmonicity of $\abs{\varphi}$ in $\X$.

Since the Jacobian $J_h=\abs{h_z}^2-\abs{h_{\bar z}}^2$
is nonnegative, it follows from~\eqref{forphi} that $h_{\bar z}$ is locally bounded, specifically
\begin{equation}\label{ho2}
\abs{h_{\bar z}(z)} \le \sqrt{\abs{\varphi(z)}} \le \frac{\norm{\varphi}^{1/2}}{\sqrt{\pi} \dist(z,\partial \X)}.
\end{equation}
The boundedness of $h_{\bar z}$ implies that $\abs{h_z}$ is locally in $BMO$, and consequently $h\in \W^{1,p}_{\rm loc}(\X)$ for every $1<p<\infty$.
A similar argument was carried out in~\cite{Mo} in a somewhat different context. However, it does not yield the Lipschitz continuity of $h$, which
we will prove by an entirely different method.

The proof of Theorem~\ref{MainTh} is preceded by several lemmas. We denote the average value of a function by an integral sign with a dash.
The normal and tangential derivatives of $h$ are defined as
\begin{equation}\label{hnt}
h_N=\frac{1}{\abs{z}}(z h_z+\bar z h_{\bar z})\qquad
h_T=\frac{i}{\abs{z}}(z h_z-\bar z h_{\bar z}).
\end{equation}

\begin{lemma}\label{avelip}
Suppose that $0\in \X$ and $h(0)=0$. Let $R=\dist(0,\partial\X)$. Then the circular mean
\begin{equation}
S(\rho):=\frac{1}{2\pi \rho} \int_{\T_\rho} h=\dashint_{\T_\rho} h, \qquad 0< \rho <R
\end{equation}
is a locally Lipschitz function of $\rho$.  Specifically we have
\begin{equation}\label{avelip1}
\abs{S'(\rho)} \le \frac{2\, \norm{\varphi}^{1/2}}{\sqrt{\pi}(R-\rho)}, \quad 0<\rho<R.
\end{equation}
As a consequence,
\begin{equation}\label{avelip1int}
\abs{S(\rho)} \le \frac{2\, \rho \norm{\varphi}^{1/2}}{\sqrt{\pi}(R-\rho)}, \quad 0<\rho<R.
\end{equation}
\end{lemma}

\begin{proof}  Note that $S$ is an absolutely continuous
function of $\rho$, even more $S\in \W_{\loc}^{1,2}(0,R)$.
Therefore, we can differentiate with respect to $\rho$ for a.e. $\rho\in (0,R)$ to obtain
\begin{equation}\label{avve1}
S'(\rho) = \dashint_{\T_\rho} h_N
= \dashint_{\T_\rho} \frac{1}{\rho} (z h_z+ \bar z h_{\bar z}).
\end{equation}
Combining~\eqref{avve1}  and the identity
\[
\dashint_{\T_\rho} \frac{1}{\rho} (z h_z - \bar z h_{\bar z})
= -i\,\dashint_{\T_\rho} h_T =0,
\]
yields
\begin{equation}\label{avelip3}
\dashint_{\T_\rho} h_N = \frac{2}{\rho}\,\dashint_{\T_\rho}  \bar z \, h_{\bar z}.
\end{equation}
This together with~\eqref{ho2} implies~\eqref{avelip1} and the integration yields~\eqref{avelip1int}.
\end{proof}

\begin{lemma}\label{equi} Suppose $0\in \X$. For almost every $\rho$, $0<\rho<\dist(0,\partial\X)$, we have
\begin{equation}\label{equi0}
\int_{\T_\rho}\abs{h_N}^2 = \int_{\T_\rho}\abs{h_T}^2 = \frac{1}{2}\int_{\T_\rho} \abs{Dh}^2.
\end{equation}
\end{lemma}

\begin{proof} Using the identities~\eqref{hnt} we find that
\begin{equation}\label{equi2}
\abs{h_N}^2+\abs{h_T}^2 = \abs{Dh}^2
\end{equation}
and
\begin{equation}\label{equi1}
\abs{h_N}^2-\abs{h_T}^2 = \frac{4}{\abs{z}^2} \re(z^2 h_z \overline{h_{\bar z}}).
\end{equation}
Integration of~\eqref{equi1} over $\T_\rho$ shows that
\begin{equation}\label{equi3}
\begin{split}
\int_{\T_\rho} (\abs{h_N}^2-\abs{h_T}^2)\,\abs{dz}
& =\frac{4}{\abs{z}^2}\re\int_{\T_\rho} z^2 \varphi(z)\,\abs{dz} \\
&=\frac{4}{\abs{z}}\im\int_{\T_\rho} z \varphi(z)\,dz=0.
\end{split}
\end{equation}
From~\eqref{equi2} and~\eqref{equi3} we obtain~\eqref{equi0}.
\end{proof}

\begin{lemma}\label{connected}
Suppose $h\in \mathfrak D (\X, \Y)$ and $w\in \partial \Y$. If the set  $h^{-1}(w)$ intersects some domain $\Omega \Subset \X$, then it also intersects $\partial \Omega$.
\end{lemma}

\begin{proof} Let $\Omega \Subset \X$ be a domain such that $h^{-1} (w) \cap \Omega$ contains a point $a$. Consider a $cd$-convergent sequence of homeomorphisms $h_j\to h$. For each $j$
\[\min_{\partial \Omega}\abs{h_j - w}\le \abs{h_j(a)-w}, \quad \mbox{because } w\not\in h_j(\Omega). \]
Since the convergence is uniform on compact sets, letting $j\to \infty$, we conclude
\[\min_{\partial \Omega}\abs{h-w}\le \abs{h(a)-w}=0. \qedhere\]
\end{proof}

Next we require a Fourier series lemma that can be viewed as a Bonnesen-type inequality~\cite{Os}.
\begin{lemma}\label{fs}
Let $f\colon \T\to\C$ be a function in $\mathscr W^{1,2}(\T)$. Expand it into the Fourier series
$f(e^{i\theta})=\sum\limits_{n\in\Z} c_n e^{i n\theta}$. If
\begin{equation}\label{fs1}
\max_{\T} \abs{f-c_0}\ge 2 \min_{\T} \abs{f-c_0},
\end{equation}
then
\begin{equation}\label{fs2}
\sum_{n\in\Z} n\abs{c_n}^2 \le \frac{99}{100} \sum_{n\in\Z} n^2\abs{c_n}^2.
\end{equation}
\end{lemma}

Indeed, inequality~\eqref{fs2}  is apart from the better factor $\nicefrac{99}{100}$,  a form of isoperimetric inequality in the plane.
We gain this better factor because of the assumption~\eqref{fs1} which can be interpreted as saying  that the image of $f$ is far from being
a circle centered at $c_0$.

\begin{proof} We may assume that $f$ is nonconstant. Normalize $f$ so that  $c_0=0$ and
$\max\limits_{\T} \abs{f-c_0}=1$. Clearly,
\begin{equation}\label{fs5}
\abs{c_1}\le 1.
\end{equation}
We need a lower bound for $\abs{c_n}$ as well. To this end, consider
\[
g(e^{i\theta})=f(e^{i\theta})-c_1 e^{i\theta}
\]
and observe that
\begin{equation}\label{fs3}
\max_{\T}\,\abs{g} \ge \max_{\T}\Abs{\abs{f}-\abs{c_1}} \ge
\frac{1}{2}\big(1-\min_{\T}\,\abs{f}\big)\ge \frac{1}{4}
\end{equation}
by virtue of~\eqref{fs1}. On the other hand,
\[
\max_{\T}\, \abs{g}^2 \le \bigg(\sum_{n\ne 1}\abs{c_n} \bigg)^2 \le
\bigg(\sum_{n\ne 1}\frac{1}{n^2}\bigg) \bigg(\sum_{n\ne 1}n^2\abs{c_n}^2 \bigg) \le 3 \sum_{n\ne 1}n^2\abs{c_n}^2
\]
which together with~\eqref{fs3} yield
\begin{equation}\label{fs4}
\sum_{n\ne 1}n^2\abs{c_n}^2 \ge \frac{1}{48}.
\end{equation}
Using~\eqref{fs5} and~\eqref{fs4} we arrive at~\eqref{fs2} as follows.
\[\begin{split}
\sum_{n\in\Z} \left(n-\frac{99}{100}n^2\right)\abs{c_n}^2
&\le \frac{1}{100}\abs{c_1}^2 + \sum_{n\ne 1}\left(\frac{1}{2}n^2-\frac{99}{100}n^2\right)\abs{c_n}^2 \\
&\le \frac{1}{100} -\frac{49}{100}\cdot \frac{1}{48}\le 0. \qedhere
\end{split}\]
\end{proof}

Let the ``good'' part of $\X$ be $G=h^{-1}(\Y)$ and the ``bad'' part be $B=\X \setminus G$. By Theorem~\ref{thharm} the restriction of $h$ to $G$ is a harmonic diffeomorphism onto $\Y$.
Thus the case $G=\X$ is trivial. From now on we assume that $B$ is nonempty. This is only possible if $\varphi\not\equiv 0$.

\begin{proposition}\label{mpro}
Suppose that $0\in \X$ and $h(0)=0\in \partial \Y$.
Then for a.e. $0<\rho<\dist(0,\partial\X)$, we have
\begin{equation}\label{mpro1}
\iint_{B_\rho} J_h \le \frac{99}{100}\, \frac{\rho}{2}\int_{\T_\rho}\abs{h_T}^2 + 4\pi \abs{c_0}^2,
\qquad \text{where } \ c_0=\dashint_{\T_\rho} h.
\end{equation}
\end{proposition}

\begin{proof} Fix $\rho$ such that the restriction of $h$ to $\T_\rho$ is in $\mathscr W^{1,2}(\T_\rho)$.
Let $M=\max\limits_{\T_\rho}\abs{h-c_0}$. The change of variables formula~\eqref{enjoyable} implies
\begin{equation}\label{mpro2}
\iint_{B_\rho} J_h \le \pi M^2.
\end{equation}
Thus we may assume $M > 2\abs{c_0}$; otherwise~\eqref{mpro1} is immediate from~\eqref{mpro2}.
By Lemma~\ref{connected} the mapping $h$ assumes the value $0$ on $\T_\rho$, hence
\[
\min_{\T_\rho}\,\abs{h-c_0}\le  \abs{c_0}<\frac{M}{2}.
\]
Thus, Lemma~\ref{fs} applies to the restriction of $h$ onto $\T_\rho$. The estimate~\eqref{fs2} reads as
\[
\iint_{B_\rho} J_h \le \frac{99}{100}\frac{\rho}{2}\int_{\T_\rho}\abs{h_T}^2
\]
which implies~\eqref{mpro1}.
\end{proof}

\begin{proof}[Proof of Theorem~\ref{MainTh}]
We will show that for every point $a\in \X$
\begin{equation}\label{ggo}
\limsup_{\rho\searrow 0} \dashiint_{B_\rho(a)} \abs{Dh}^2 \le \frac{16000}{\pi \dist^2(a,\partial \X)}\iint_\X \abs{Dh}^2.
\end{equation}
We may assume $a=0=h(a)$. Let $R=\dist(0,\partial\X)$.

\textbf{Case 1.} $0\in \partial \Y$. The first step is to rewrite~\eqref{mpro1} as a differential inequality for the function
\[
\EE(\rho): = \iint_{B_\rho}\abs{Dh}^2,\qquad 0< \rho< \frac{R}{2}.
\]
Here we restrict ourselves to $0<\rho<\nicefrac{R}{2}$ which yields $\dist (z, \partial \X) \ge \nicefrac{R}{2}$ for $z\in B_\rho$.
Since $\abs{Dh}^2=2J_h+4\abs{h_{\bar z}}^2$, applying the inequality~\eqref{ho2} we obtain
\begin{equation}\label{o18}
\EE(\rho) \le 2\iint_{B_\rho} J_h +  \frac{4\rho^2\norm{\varphi}}{(\nicefrac{R}{2})^2}
= 2\iint_{B_\rho} J_h +  \frac{16\rho^2\norm{\varphi}}{R^2}.
\end{equation}
Next, the integral on the right is estimated using~\eqref{mpro1},~\eqref{avelip1int} and Lemma~\ref{equi}:
\begin{equation}\label{o19}
2\iint_{B_\rho} J_h \le \frac{99}{100}\rho\int_{\T_\rho}\abs{h_T}^2 + 8\pi \abs{c_0}^2 \le
\frac{99}{200}\rho\int_{\T_\rho}\abs{Dh}^2 + \frac{128\rho^2 \norm{\varphi} }{R^2}.
\end{equation}
Combining~\eqref{o18} and~\eqref{o19} we obtain
\begin{equation}\label{o20}
\EE(\rho) \le \frac{99}{200}\rho\, \EE'(\rho) + \frac{144 \rho^2 \norm{\varphi}}{R^2} \quad \mbox{  for a.e. }0<\rho<\frac{R}{2}.
\end{equation}

For notational simplicity we introduce the constant $q=\nicefrac{200}{99}$. Inequality~\eqref{o20} yields
\begin{equation}\label{o21}
\frac{d}{d\rho} \left(\rho^{-q}\EE(\rho)\right)
=\frac{\rho \,  \EE'(\rho)-q\, \EE(\rho)}{\rho^{\, q+1}}
\ge -\frac{144 \, q\,  \norm{\varphi}}{R^2 \rho^{\, q-1}}.
\end{equation}
Integrate~\eqref{o21} over the interval $(\rho,\nicefrac{R}{2})$ to obtain
\begin{equation}\label{o22}
\frac{1}{\rho^{\, q}}\EE(\rho) \le \frac{2^q}{R^q}\EE (\nicefrac{R}{2}) + \frac{144 \,q\, \norm{\varphi}}{(q-2)R^2 \rho^{\, q-2}}
= \frac{2^q}{R^{\, q}}\EE (\nicefrac{R}{2}) + \frac{14400 \, \norm{\varphi}}{R^2 \rho^{\, q-2}}.
\end{equation}
Finally, multiply~\eqref{o22} by $\rho^{q-2}$ and rewrite it as
\begin{equation}\label{pro2}
\dashiint_{B_\rho} \abs{Dh}^2 \le \frac{2^q }{\pi R^2} \left(\frac{\rho}{R}\right)^{q-2} \iint_{B_{R/2}} \abs{Dh}^2 + \frac{14400 \, \norm{\varphi}}{\pi R^2}.
\end{equation}
We further simplify~\eqref{pro2} using the pointwise inequality $4\abs{\varphi}\le \abs{Dh}^2$.
\begin{equation}\label{pro2b}
\dashiint_{B_\rho} \abs{Dh}^2 \le \left(2^q \left(\frac{\rho}{R}\right)^{q-2}+3600\right)\frac{1}{\pi R^2} \iint_\X \abs{Dh}^2.
\end{equation}
This yields~\eqref{ggo} even with a better constant.

\textbf{Case 2.} $0\in \Y$; that is, $0\in G=h^{-1}(\Y)$.
Let $r=\dist(0,\partial G)\le \dist (0,\partial \X)= R$. Since $h$ is harmonic in $G$, the subharmonicity
of $\abs{Dh}^2$ yields
\begin{equation}\label{pro3}
\abs{Dh(0)}^2\le \dashiint_{B_r}\abs{Dh}^2 \le \frac{1}{\pi r^2}\iint_{\X}\abs{Dh}^2.
\end{equation}
If $R<60\, r$, then~\eqref{pro3} already implies~\eqref{ggo}.
Otherwise pick $\zeta\in \X\setminus G$ such that $\abs{\zeta}=r$. Clearly
\begin{equation}\label{o25}
\dashiint_{B_r}\abs{Dh}^2 \le 4\,\dashiint_{B_{2r}(\zeta)}\abs{Dh}^2
\end{equation}
and the righthand side of~\eqref{o25} can be estimated by applying~\eqref{pro2b} to $B_{2r}(\zeta)$.
Since $\dist(\zeta,\partial \X)\ge R-r \ge 59\, r$, we may apply inequality~\eqref{pro2b} with $\rho=2r$ and $R-r$ in place of $R$. This gives the estimate
\begin{equation}\label{o26}
\begin{split}
\dashiint_{B_{2r}(\zeta)}\abs{Dh}^2 &\le \left(2^q \left(\nicefrac{2}{59}\right)^{q-2}+3600\right)\frac{(\nicefrac{60}{59})^2}{\pi R^2} \iint_\X \abs{Dh}^2 \\
&\le \frac{4000}{\pi R^2} \iint_\X \abs{Dh}^2.
\end{split}
\end{equation}
Combining~\eqref{pro3},~\eqref{o25}, and~\eqref{o26} gives the inequality~\eqref{ggo}.

This estimate yields the pointwise inequality~\eqref{addstar} at the Lebesgue points of $\abs{Dh}^2$, completing the proof of Theorem~\ref{MainTh}.
\end{proof}

\section{$\mathscr C^1$-smoothness of minimal deformations, proof of Theorem~\ref{smooth}}\label{secsmo}


Throughout this section the following standing assumptions are made on the mappings under considerations:
$\X$ and $\Y$ are $\ell$-connected bounded domains  and $h \colon \X \to \overline{\Y}$ is a $\mathscr C^1$-smooth Hopf deformation, $h\in \mathfrak D (\X, \Y)$; that is,
\[h_z \overline{h_{\bar z}} = \varphi \not\equiv 0 \qquad \mbox{in } \X\]
where $\varphi$ is a holomorphic function. Recall that the \emph{convex part} of the boundary of a domain $\Y$ is
\begin{equation}\label{eqconv}
\partial_c\Y=\{w\in \partial \Y\colon B_r(w)\cap \Y \quad \text{is convex for some }\ r>0\}.
\end{equation}
 We designate the domain of regular points of the quadratic differential $\varphi (z)\, \dtext z^2$ by
\[\X_\circ = \X \setminus \{\mbox{zeros of } \varphi\}.\]
Additional assumptions on $h$ will be explicitly stated when needed. The proof of Theorem~\ref{smooth} proceeds by a number of claims and Lemma~\ref{bes}.

{\bf Claim 1.} {\it Let $0\in h(\X)$, then each horizontal arc $\alpha \subset \overline{\alpha} \subset \X_\circ$ contains at most a finite number of zeros of $h$.}
\begin{proof}
By virtue of~\eqref{hHhVJh} we have a lower bound for the horizontal derivative
\begin{equation}\label{lbhd}
\abs{ \partial_{_\mathsf H} h }\ge 2 \sqrt{\abs{\varphi}}>0 \qquad \mbox{for all } z \in \overline{\alpha}.
\end{equation}
If $h\colon \overline{\alpha} \to \C$ had an infinite number of zeros there would be an accumulation point of zeros. Since $h$ is $\mathscr C^1$-smooth along $\overline{\alpha}$, its horizontal derivative would vanish at the accumulation point of zeros, in contradiction with~\eqref{lbhd}.
\end{proof}

{\bf Claim 2.} {\it Suppose $0\in \partial \Y$. Then the preimage $h^{-1}(0)$ is covered by a countable number of vertical trajectories and critical points.}
\begin{proof}
The domain $\X_\circ$ can be expressed as countable union of $\varphi$-rectangles compactly contained in $\X_\circ$, see \S\ref{secqd} for the definition and consideration of $\varphi$-rectangles. Let $\mathfrak R$ be one such $\varphi$-rectangle. It suffices to show that the set $h^{-1}(0) \cap \mathfrak R$ can be covered by a countable number of vertical trajectories. Let $\alpha$ and $\beta$ denote the horizontal edges of $\mathfrak R$. These are horizontal arcs of $\varphi \, \dtext z^2$ compactly contained in $\X_\circ$. In view of Claim 1, we may and do assume that $h$ does not vanish on the horizontal edges of $\mathfrak R$, for otherwise we can replace $\mathfrak R$ by a finite number of parallel subrectangles of $\mathfrak R$; simply cut $\mathfrak R$ along vertical trajectories crossing the horizontal edges $\alpha$ and $\beta$ at  the zeros of $h$.  With this assumption the proof of Claim 2 will be completed by showing that inside of $\mathfrak R$ there are no zeros of $h$. Indeed, suppose otherwise. Then, by Lemma~\ref{connected}, the set $h^{-1}(0)$ intersects at least one of two vertical edges of the rectangle $\mathfrak R$. By the same reasoning $h^{-1}(0)$ intersects one of two vertical edges of any subrectangle obtained by continuously compressing $\mathfrak R$ with vertical trajectories. This procedure results in a family of vertical arcs intersecting $h^{-1}(0)$ whose union has positive measure, in contradiction with Lemma~\ref{lily}.
\end{proof}
It should be mentioned that a result similar to Claim 2 was proved in~\cite{Ku}  for energy-minimizing mappings between closed surfaces that are smooth except for isolated singularities.

Recall that Theorem~\ref{smooth} deals with minimal deformations. These are special mappings for which the Hopf differential $h_z \overline{h_{\bar z}}\, \dtext z^2$ is real along $\partial \X$, provided $\partial \X$ is $\mathscr C^1$-smooth. It involves no loss of generality in the proof of Theorem~\ref{smooth} to assume that $\partial \X$ is indeed $\mathscr C^1$-smooth. For, if necessary, we could transform $\X$ conformally onto, say, a circular Schottky domain.
Our next lemma tells us that no point in $\X$ can be mapped into a convex part $\partial_c\Y$ of $\partial \Y$, see~\eqref{eqconv} and Figure~\ref{fig4}.

\begin{lemma}\label{bes}
Suppose $\partial \X$ is $\mathscr C^1$-smooth and $\varphi (z)\, \dtext z^2$ is real on the boundary of $\X$, see Definition~\ref{realqd} in \S\ref{secqd}. Assume that $h$ is $\mathscr C^1$-smooth in $\X$. Then
\begin{equation}\label{8dstar}
h(\X) \cap \partial_c \Y = \varnothing.
\end{equation}
\end{lemma}
\begin{proof}
Assume, to the contrary, that the set
\begin{equation}
B_c:=h^{-1} (\partial_c \Y) \ne \varnothing.
\end{equation}
Consider the ``good'' subdomain $G= h^{-1}(\Y) \subset \X$. It is indeed connected by Theorem~\ref{thharm}. Let $E$ be the set of points $z\in \partial G \cap \X$ for which there is a disk $D=D_z \subset G$ such that $z\in \partial D$. Note that $E$ is dense in $\partial G \cap \X$. Our first step is to prove that

{\bf Step 1.} {\it We have
\begin{equation}\label{fstep}
h(E)\cap \partial_c\Y = \varnothing.
\end{equation}
Hence, more generally,}
\begin{equation}\label{eq84}
h(\partial G \cap \X) \cap \partial_c \Y = \varnothing.
\end{equation}
\begin{proof}
Suppose to the contrary that $w = h(z)\in \partial_c\Y$ for some $z\in \X$. Since $\Y$ is convex at $w$
the function $u:=\re (e^{i\theta} h)$, for some constant $0\le \theta \le 2 \pi$, attains a local maximum at $z$, which yields
\begin{equation}\label{you}
 \nabla u(z)=0.
\end{equation}
On the other hand, there is a disk $D\subset G$ such that $z\in\partial D$, $h(D) \subset h(G) \subset \Y$ and  $u(z)>u(\zeta)$ for all $\zeta \in D$. The  Hopf boundary point lemma (\cite[Lemma 3.4]{GTb}, \cite{HopfE})
implies that the inner normal derivative of $u$ at $z\in\partial D$ is strictly negative. This contradicts~\eqref{you}. Finally, since $\partial_c \Y$ is open in $\partial \Y$ and $E$ is dense in $\partial G \cap \X$, the claim~\eqref{eq84} follows.
\end{proof}

{\bf Step 2.} {\it The set $B_c= h^{-1} (\partial_c \Y) \subset \X$ is open.}
\begin{proof}
By virtue of~\eqref{eq84} the set $B_c$ is contained in the interior of $\X \setminus h^{-1} (\Y)$. Since $\partial_c \Y$ is open in $\partial \Y$, it follows that $B_c$ is open in $\X$.
\end{proof}

{\bf Step 3.} {\it Every noncritical vertical trajectory $\Gamma$ that intersects $B_c$ is a connected component of $h^{-1} (w)$ for some $w\in \partial_c \Y$.}
\begin{proof}
The intersection of $\Gamma$ with $B_c$ is open in $\Gamma$. On the other hand the vertical derivative  $\partial_{_\mathsf V} h$ vanishes a.e. on  $\X\setminus h^{-1} (\Y)$. Indeed,
by~\cite[Lemma~3.8]{IKKO} $J_h$ vanishes a.e. on $\X\setminus h^{-1} (\Y)$
and by~\eqref{becareful} the vertical derivative satisfies $\abs{\partial_{_\mathsf V} h} \le \sqrt{\abs{J_h}}=0$.
Therefore, $h$ is constant along $\Gamma$, say equal to $w\in \partial_c \Y$. On the other hand, by Claim 2, the preimage $h^{-1}(w)$ is covered by a countable union of vertical trajectories and critical points. We find that $\Gamma$ is indeed one of the connected components of $h^{-1}(w)$.
\end{proof}

{\bf Step 4.} {\it A noncritical vertical trajectory $\Gamma$ that intersects $B_c$ approaches two different boundary components of $\X$.}
\begin{proof}
Recall from the introduction, \S\ref{sec101}, the components $\mathfrak X_1, \dots, \mathfrak X_\ell$ of $\mathfrak X = \partial \X$. Suppose otherwise, that both ends of $\Gamma$ approach the same boundary component, say $\mathfrak X_1$. By Lemma~\ref{callsomething} the vertical trajectory $\Gamma$ separates
two other  boundary components of $\X$, say $\mathfrak X_2$ and $\mathfrak X_3$. Let $C\subset \Y$ be a continuum  connecting  the components $\Upsilon_2$ and $\Upsilon_3$ of $\partial \Y$, in accordance with the correspondence in~\eqref{Corr}.  Then $h^{-1}(C)$ is a continuum~\cite[Lemma 3.7]{IKKO} that connects $\mathfrak X_2$ and $\mathfrak X_3$ in $\X\setminus \Gamma$, a contradiction.
\end{proof}
Now we have the required contradiction that proves Lemma~\ref{bes}. Indeed $\Gamma$, being a connected component of $h^{-1}(w) \in \partial \Y$, cannot approach different components of $\partial \X$; otherwise its image $h(\Gamma)$ would approach different components of $\partial \Y$
by virtue of~\eqref{Corr}.
\end{proof}

To complete the proof of Theorem~\ref{smooth} we need one more lemma.
\begin{lemma}\label{bes2}
If $\X \setminus h^{-1} (\Y)$ has $\sigma$-finite $1$-dimensional Hausdorff measure, then  the set  \mbox{$\X \setminus h^{-1} (\Y)$}  is empty and $h$ is a diffeomorphism of $\X$ onto $\Y$.
\end{lemma}
\begin{proof} The function $h_z$ is continuous in $\X$ and holomorphic in $h^{-1}(\Y)$.
By a theorem of Besicovitch~\cite[Theorem 2]{Be}  sets of  $\sigma$-finite $1$-dimensional Hausdorff measure are removable for holomorphic continuous functions. Therefore, $h_z$ is holomorphic in $\X$. It follows that $h$ is harmonic in $\X$, which by Proposition~\ref{harmdiff} implies that $h$ is a diffeomorphism.
\end{proof}

\begin{proof}[Proof of Theorem~\ref{smooth}.]
Combining Lemmas~\ref{bes} and~\ref{bes2} we find that the set $\X\setminus h^{-1}(\Y)$ has $\sigma$-finite
1-dimensional measure and therefore $h$ is a diffeomorphism by Lemma~\ref{bes2}.
\end{proof}

\section{Open questions}

The Hopf-deformations in Example~\ref{exsqueeze} are slightly better than Lipschitz regular. Namely, their energy functions are continuous. This raises the following question.

\begin{question}\label{ask0}
Suppose that $h$ is an  energy-minimal deformation. Is $|Dh|$ continuous?
\end{question}

The reader may notice that the domain $\Y$ in Theorem~\ref{smooth} cannot be $\mathscr C^1$-smooth. Also, the butterfly domain $\Y$ of Example~\ref{exsqueeze} has nonsmooth boundary.
These observations lead to the following question.

\begin{question}\label{ask}
Are energy-minimal deformations $\mathscr C^1$-smooth when the boundary of $\Y$ is smooth?
\end{question}

In Example~\ref{exhammering} the energy-minimal deformation $h$ belongs to $\mathscr C^{1,1}(\mathbb X)$ but not to $\mathscr C^2(\mathbb X)$. Thus one may expect the $\mathscr C^{1,1}$-regularity to hold in Question~\ref{ask}.

\bibliographystyle{amsplain}

\begin{thebibliography}{99}

\bibitem{Ah}
L. V. Ahlfors,
\textit{Quasiconformal deformations and mappings in $\mathbb{R}^{n}$},
J. Analyse Math. \textbf{30} (1976), 74--97.


\bibitem{AIMb}
K. Astala, T. Iwaniec, and G. Martin, \textit{Elliptic partial differential equations and quasiconformal mappings in the plane},  Princeton University Press, Princeton, NJ, 2009.


\bibitem{AIM}
K. Astala, T. Iwaniec, and G.  Martin,  \textit{Deformations of annuli with
smallest mean distortion}. Arch. Ration. Mech. Anal. {\bf 195}, no.~3 (2010), 899--921.


\bibitem{Be}
A. S. Besicovitch, \textit{On sufficient conditions for a function to be analytic, and on behaviour of analytic functions in the neighbourhood of non-isolated singular
points}, Proc. London Math. Soc. (2) \textbf{32} (1931), no.~1, 1--9.



\bibitem{Ch}
J. Chen,
\textit{On energy minimizing mappings between and into singular spaces},
Duke Math. J.  \textbf{79}  (1995),  no.~1, 77--99.

\bibitem{Ch1}
Y. W. Chen,  \textit{Discontinuity and representations of minimal surface solutions}, Proceedings of the conference on differential equations (dedicated to A. Weinstein), pp. 115--138. University of Maryland, College Park, MD (1956).

\bibitem{CL}
E. F. Collingwood and A. J. Lohwater,
\textit{The theory of cluster sets}, Cambridge University Press, Cambridge, 1966.

\bibitem{Cob}
R. Courant, \textit{Dirichlet's principle, conformal mapping, and minimal surfaces}, 
With an appendix by M. Schiffer.  Springer-Verlag, New York-Heidelberg, 1950.

\bibitem{DM}
G. Daskalopoulos and Ch. Mese,
\textit{Harmonic maps between singular spaces. I},
Comm. Anal. Geom. \textbf{18}  (2010),  no.~2, 257--337.

\bibitem{DHTb}
U. Dierkes,  S. Hildebrandt,  and A. J.  Tromba,  \textit{Global Analysis of Minimal Surfaces}, Springer-Verlag, Berlin, 2010.

\bibitem{EL2}
J. Eells and L. Lemaire, \textit{Selected topics in harmonic maps}, CBMS Regional Conference Series in Mathematics, 50.  American Mathematical Society, Providence, RI, 1983.

\bibitem{EL3}
J. Eells and L. Lemaire, \textit{Another report on harmonic maps}, Bull. London Math. Soc. \textbf{20} (1988), no.~5, 385--524.

\bibitem{GMb}
J. B. Garnett and D. E. Marshall,
\textit{Harmonic measure}, Cambridge Univ. Press, Cambridge, 2005.

\bibitem{GTb}
D. Gilbarg and N. S. Trudinger,
\textit{Elliptic partial differential equations of second order}. Second edition.
Springer-Verlag, Berlin, 1983.

\bibitem{GS}
M. Gromov and R. Schoen,
\textit{Harmonic maps into singular spaces and $p$-adic superrigidity for lattices in groups of rank one},
Inst. Hautes \'Etudes Sci. Publ. Math. No. \textbf{76} (1992), 165--246.

\bibitem{He}
F. H\'elein,
\textit{Hom\'eomorphismes quasi conformes entre surfaces riemanniennes},
C. R. Acad. Sci. Paris S\'er. I Math.  \textbf{307}  (1988),  no.~13, 725--730.

\bibitem{Heb}
F. H\'elein,
\textit{Harmonic maps, conservation laws and moving frames}, 2nd edition. Cambridge University Press, Cambridge, 2002.



\bibitem{HopfE}
E. Hopf,
\textit{A remark on linear elliptic differential equations of second order},
Proc. Amer. Math. Soc. \textbf{3} (1952), 791--793.

\bibitem{HopfH}
H. Hopf,
\textit{Differential geometry in the large},
Notes taken by Peter Lax and John Gray. With a preface by S. S. Chern. Lecture Notes in Mathematics, 1000. Springer-Verlag, Berlin, 1983.

\bibitem{IKKO}
T. Iwaniec, N.-T. Koh, L. V. Kovalev, and J. Onninen,
\textit{Existence of energy-minimal diffeomorphisms between doubly connected domains}, Invent. Math. \textbf{186} (2011), no. 3, 667--707.

\bibitem{IKOhopf}
T. Iwaniec, L. V. Kovalev, and J. Onninen, \textit{Hopf differentials and smoothing Sobolev homeomorphisms},  {Int. Math. Res. Not. IMRN \textbf{2012} (2012), no. 14, 3256--3277.} 


\bibitem{IMb}
T.  Iwaniec and G. Martin, \textit{Geometric function theory and non-linear analysis},  Oxford University Press, New York, 2001.

\bibitem{IOt}
T. Iwaniec and J. Onninen, \textit{Deformations of finite conformal energy: boundary behavior and limit theorems},  Trans. Amer. Math. Soc. 363 (2011), no. 11, 5605--5648.

\bibitem{IS}
T. Iwaniec \and V. \v Sver\'ak,
\textit{On mappings with integrable dilatation},
Proc. Amer. Math. Soc. \textbf{118} (1993), no.~1, 181--188.

\bibitem{Jono}
J. Jost, \textit{A note on harmonic maps between surfaces},
Ann. Inst. H. Poincar\'e Anal. Non Lin\'eaire 2 (1985), no. 6, 397--405.

\bibitem{Job}
J. Jost, \textit{Two-dimensional geometric variational problems}, John Wiley \& Sons, Ltd., Chichester, 1991.

\bibitem{Jo96}
J. Jost,
\textit{Generalized harmonic maps between metric spaces}, in ``Geometric analysis and the calculus of variations'', 143--174, Int. Press, Cambridge, MA, 1996.

\bibitem{Jo}
J. Jost, \textit{Generalized Dirichlet forms and harmonic maps},
Calc. Var. Partial Differential Equations \textbf{5} (1997), no.~1, 1--19.

\bibitem{KS}
N. J. Korevaar and R. M. Schoen,
\textit{Sobolev spaces and harmonic maps for metric space targets},
Comm. Anal. Geom.  \textbf{1}  (1993),  no.~3--4, 561--659.

\bibitem{Ku}
E. Kuwert, \textit{Harmonic maps between flat surfaces with conical singularities}, Math.~Z. \textbf{221} (1996), no.~3, 421--436.

\bibitem{LWb}
F. Lin and Ch. Wang, \textit{The analysis of harmonic maps and their heat flows}, World Scientific Publishing Co. Pte. Ltd., Hackensack, NJ, 2008.

\bibitem{MM}
V. Markovi\'c and M. Mateljevi\'c,
\textit{A new version of the main inequality and the uniqueness of harmonic maps},
J. Anal. Math. \textbf{79} (1999), 315--334.

\bibitem{Mc}
L. F.  McAuley,   \textit{Some fundamental theorems and problems related to monotone mappings},  Proc. First Conf. on Monotone Mappings and Open Mappings, SUNY at Binghamton, Binghamton, N.Y. (1971), 1--36.


\bibitem{Mo}
R. Moser,
\textit{On a variational problem with non-differentiable constraints},
Calc. Var. Partial Differential Equations \textbf{29} (2007), no.~1, 119--140.


\bibitem{Os}
R. Osserman, \textit{Bonnesen-style isoperimetric inequalities}. Amer. Math. Monthly {\bf 86} (1979), no. 1, 1--29.


\bibitem{Rab}
T. Rad\'o,  \textit{Length and Area},  Amer. Math. Soc., 1948.

\bibitem{Ranb}
T. Ransford, \textit{Potential theory in the complex plane},
Cambridge University Press, Cambridge, 1995.




\bibitem{RS}
E. Reich and K. Strebel,
\textit{On quasiconformal mappings which keep the boundary points fixed},
Trans. Amer. Math. Soc. \textbf{138} (1969), 211--222.




\bibitem{SS}
J. Sivaloganathana and S. J. Spector, \textit{On irregular weak solutions of the energyÐmomentum equations},  
Proc.  R. Soc.  Edinb. A \textbf{141} (2011), 193--204.

\bibitem{Stb}
K. Strebel, \textit{Quadratic differentials}, Springer-Verlag, Berlin, 1984.

\bibitem{Ta}
A. Taheri,  \textit{Quasiconvexity and uniqueness of stationary points in the multi-dimensional calculus of variations}, Proc. Amer. Math. Soc. \textbf{131} (2003), no.~10, 3101--3107. 

\bibitem{Wa}
Ch. Wang,
\textit{Energy minimizing maps to piecewise uniformly regular Lipschitz manifolds},
Comm. Anal. Geom. \textbf{9} (2001), no.~4, 657--682.

\bibitem{Wh}
G. T. Whyburn, \textit{Analytic Topology},
 American Mathematical Society, New York, 1942.

\bibitem{Yo}
J. W. T. Youngs, \textit{The topological theory of Fr\'echet surfaces},
Ann. of Math. (2) \textbf{45} (1944), 753--785.

\end{thebibliography}

\end{document}